\documentclass[preprint,review,12pt]{elsarticle}

\usepackage{amssymb}
\usepackage{subcaption}
\usepackage{graphicx}
\usepackage{verbatim}
\usepackage{mathtools}
\usepackage{hyperref}
\usepackage{mathrsfs}
\usepackage{enumitem}
\usepackage{color}
\usepackage{booktabs}
\usepackage{natbib}
\usepackage[algoruled,boxed,lined]{algorithm2e}
\usepackage{algorithmic}
\usepackage[normalem]{ulem}
\usepackage{tikz}
\usepackage{floatrow}
\DeclareRobustCommand\sampleline[1]{%
   \tikz\draw[#1][line width=0.5mm] (0,0) (0,\the\dimexpr\fontdimen22\textfont2\relax)
   -- (2em,\the\dimexpr\fontdimen22\textfont2\relax);%
    }

\def\bf#1{\boldsymbol{#1}}

\def\red#1{\textcolor{red}{#1}}
\def\blue#1{\textcolor{blue}{#1}}

\journal{Applied Energy}

\begin{document}
	
	\begin{frontmatter}
		
		\title{Particle swarm optimization of a wind farm layout with active control of turbine yaws}
		
		\author{Jeonghwan Song \fnref{cont}}
		\author{Taewan Kim \fnref{cont}}
		\author{Donghyun You \corref{cor1}}
		\address{Department of Mechanical Engineering, Pohang University of Science and Technology, 77 Cheongam-Ro, Nam-Gu, Pohang, Gyeongbuk 37673, South Korea \vspace{-0.4in}}
		
		\fntext[cont]{Authors contributed equally.}
		\cortext[cor1]{Corresponding author.}
		
		\ead{dhyou@postech.ac.kr}
		
		\begin{abstract}
		Active yaw control (AYC) of wind turbines has been widely applied to increase the annual energy production (AEP) of a wind farm. AYC efficiency depends on the wind direction and the wind farm layout because an AYC method utilizes wake deflection by yawing wind turbines. Conventional optimization of a wind farm layout assumed that the swept areas of all wind turbines are aligned perpendicular to the wind direction, thereby allowing non-optimal utilization of an AYC method. Higher AEP can be obtained by joint optimization which considers an AYC method in the layout design stage. Joint optimization of the farm layout and AYC has been difficult due to the non-convexity of the problem and the computational inefficiency. In the present study, a particle swarm optimization based method is developed for joint optimization. The layout is optimized with simultaneous consideration for yaw angles for all wind velocities to obtain a globally optimal layout. A number of random initial particles consisting of the layout and yaw angles of wind turbines reduce the initial layout dependency on the optimized layout. To deal with the challenge of large-scale optimization, the adaptive granularity learning distributed particle swarm optimization algorithm is implemented. The improvement in AEP when using a jointly optimized layout compared to a conventionally optimized layout in a real wind farm is demonstrated using the present method.

		\end{abstract}
		
		
	\end{frontmatter}

	\section{Introduction}\label{sec:Introduction}
	The annual energy production (AEP) of a wind farm is significantly reduced by wake interactions~\cite{barthelmie2009, cleve2009, barthelmie2010}. Once a wind turbine converts wind energy into electrical energy, a wake region with a lower wind speed than its surroundings is created. The power production of a wind turbine in the wake region is lower than that of an upstream wind turbine due to the velocity deficit. Therefore, the way to increase AEP of a wind farm is to reduce the overlap of the wind turbine swept areas and the wake regions as much as possible. For this purpose, two main strategies, optimization of a wind farm layout~\cite{parada2017, li2017,  mirhassani2017, wu2021, zilong2022, reddy2020, antonini2018, kuo2016, yang2019, long2020} and active wind turbine control~\cite{boulamatsis2019, cooperman2015, menezes2018, tang2018, lin2020, zong2021, song2022, song2018, ma2021} have been utilized. Layout optimization shifts the positions of wind turbines to separate downstream wind turbines from wake regions. Active wind turbine control changes the blade pitch angle and the rotor tilt and yaw angles to weaken or deflect the wake of a wind turbine. Among current active control options, active yaw control (AYC) is the most promising strategy because of its large wake deflections and the small effect on the structural loads~\cite{zong2021, dou2020, ma2021, song2022, song2018}.

    An optimal layout of a wind farm depends on the condition of the wind farm, such as the wind turbine model, the number of turbines, the land size, and the wind rose. Layout optimization often becomes a highly non-convex and complex problem, therefore meta-heuristic optimization algorithms have been used to solve. The most commonly used algorithms are evolutionary computation algorithms, including genetic algorithms and particle swarm optimization (PSO) algorithms~\cite{wu2020, wu2021, ju2019, masoudi2022, brogna2020, kirchner2018, tao2020,chowdhury2012, hou2017}. However, as the number of turbines of a layout optimization problem increases, the problem becomes non-convex and leads to non-optimal local minima. To search for a globally optimal solution to a large-scale problem of layout optimization, improved algorithms have been developed~\cite{dou2020, patel2017, pookpunt2016}.

	AYC steers the yaw angles of wind turbines to separate downstream wind turbines from wake regions. The efficiency of AYC for a wind farm depends on the wind farm layout and the wind direction. AYC is more efficient for particular wind directions where turbines are more affected by the wake. In the Horns-Rev wind farm, the highest AYC efficiency can be achieved when the wind direction is aligned with turbine rows~\cite{zong2021, dou2020}.
	
	The effect of AYC on an optimal layout has recently been studied. The conventional layout optimization problem assumed swept areas of all wind turbines are aligned perpendicular to the wind direction, which is called the greedy control strategy~\cite{parada2017, li2017,  mirhassani2017, wu2021, zilong2022, reddy2020, antonini2018, kuo2016, yang2019, long2020}. AYC can be optimized after the wind farm layout is optimized, which is defined as sequential optimization. Optimal AYC can generate power more or equal to the greedy strategy~\cite{zong2021, dou2020, ma2021}. When AYC is considered at the layout design stage, which is defined as joint optimization, the optimized layout can generate power more or equal to the sequentially optimized layout~\cite{chen2022}. However, the dimension of the joint optimization problem is significantly increased because the yaw angles of wind turbines have to be optimized for each wind velocity. Chen \textit{et al.}~\cite{chen2022} designed a decomposition-based hybrid method (DBHM) to reduce the computational costs of the joint optimization problem. Starting from the sequentially optimized layout, the optimization problem was decomposed into subproblems that optimized the layout and yaw angles for each wind velocity. The layout of the whole problem was updated based on the alternating direction method of multipliers (ADMM) using optimal solutions to subproblems. This approach successfully lowered the computational cost, but the optimal solution of DBHM highly depends on the initial layout. Also, to converge to a globally optimal solution with high accuracy, ADMM can be very slow~\cite{boyd2011}.
	
	AEP of the optimized wind farm is directly linked to investment gains, therefore finding a globally optimal layout is very important for economic reasons. The yaw angles for each wind velocity can be separately optimized, but the layout cannot be separated for each wind velocity. A decomposition-based method showed the low performance when dealing with partially separable problems~\cite{wang2020}. For partially separable and large-scale problems, such as layout optimization, algorithms that deal with the whole problem at once perform better~\cite{cheng2015, cheng2014, tan2006, liu2018, zhan2016, wang2019, wang2021}.
	
	In the present study, a PSO-based layout optimization method with simultaneous consideration of yaw angles for all wind velocities is developed to find a globally optimal layout. Since the present method does not start from one initial layout but uses a number of random initial particles consisting of layouts and yaw angles, the solution of the present method does not depend on one initial layout. Particles are updated towards the layout that simultaneously maximizes AEP for all wind velocities using information of the optimal particle of each iteration.
	
	In the present work, an adaptive granularity learning distributed particle swarm optimization (AGLDPSO) algorithm is employed, which treats a large-scale problem as a whole to avoid local minima~\cite{wang2020}. The present method shows outstanding optimization performance in solving a joint optimization problem. The optimal layout of the present method and that of the previous method in a joint optimization problem are compared. Also, in order for joint optimization to be useful for industrial sites, the joint optimization performance of the following wind farm design conditions is analyzed: 
	
\bigskip
	1. Power curve gradients of wind turbines,
	
	2. Number of wind turbines per area,
	
	3. Wind uniformity.
\bigskip

    \noindent Because a yawed wake is a three-dimensional phenomenon, an in-house code is developed to calculate three-dimensional wind farm power, and the method to solve the joint optimization is developed using the code.

	The paper is organized as follows: in Section \ref{sec:Background}, the theoretical components in the wind farm power calculation and the optimization method are explained in detail. In Section \ref{sec:Problem setup}, the configurations of the layout optimization problem are defined. The results and discussion of the problems considered are presented in Section \ref{sec:Optimization result}, followed by concluding remarks in Section \ref{sec:Concluding}.
	
	\section{Background}\label{sec:Background}
	\subsection{Yawed Gaussian wake model}\label{subsec:Yawed Gaussian wake model}
	
	The velocity of wind passing through each wind turbine needs to be calculated to obtain AEP of a wind farm. The velocity field in a wind farm depends on wake effects. Wake models which can analytically derive the wake velocity are being developed because their lower computational costs benefit layout optimization~\cite{jensen1983, frandsen2006, bastankhah2014}. To optimize AYC, a yawed wake model that considers wake steering by yaw control is needed. Dou \textit{et al.}~\cite{dou2020} developed a Gaussian-based yawed wake model which estimates the shape of the wake velocity distribution in the streamwise direction. The model equation can be written as follows:
	\begin{equation} \label{eq:wakemodelDou1}
    \begin{gathered}
    \frac{\Delta u}{u_{\infty}} = \left(1-\sqrt{1-\frac{C_{T*}\cos{\gamma}}{8\sigma_{\rm{yaw}}\sigma_z}}\right)\\
    \times\exp{\left[-\frac{1}{2\sigma_{\rm{yaw}}^2}\left(\frac{y-Y_{\rm{offset,z}}}{D\cos{\gamma}}\right)^2-\frac{1}{2\sigma_z^2}\left(\frac{z-z_{\rm{hub}}}{D}\right)^2\right]},
    \end{gathered}
	\end{equation}
		
	\noindent where 
	{
	
	\centering{$\delta_{*} = 0.607, \zeta = 0.75, k = 0.0125$}
	
	\centering{$\delta = \delta_{*}C_{T}, C_{T*} = C_{T}\cos^2{\gamma},$}
	
	\centering{$\beta = \left(1+\sqrt{1-C_{T*} \cos{\gamma}}\right)/\left(2\sqrt{1-C_{T*} \cos{\gamma}}\right),$}
	
	\centering{$\sigma_{\rm{yaw}} = kx/\left(D \cos{\gamma}\right) + \sqrt{\beta}/5,$}
	
	\centering{$\sigma_{z} = kx/D + \sqrt{\beta}/5,$}
		
	\centering{$Y_{\rm{offset}}/D = \delta\left(C_{T}\sin{\gamma}\right)^\zeta \cos^{2\zeta}{\gamma}\sqrt{x/D}+d_{rt}\sin{\gamma}/D,$}
		
	\centering{$Y_{\rm{offset,z}} = \left(Y_{\rm{offset}} - d_{rt}\sin{\gamma}\right)exp{\left[-0.5\left(z-z_{\rm{hub}}\right)^2/\left(D\sigma_z\right)^2\right]}+d_{rt}.$}
	
	}
    \noindent $D$ is the diameter of the wind turbine rotor, $z_{\rm{hub}}$ is the height of the wind turbine hub, $d_{rt}$ is the distance between the rotor center and the hub center, $\gamma$ is the yaw angle of the wind turbine, and $C_T$ is the thrust coefficient of the wind turbine when $\gamma=0$.
    
\begin{figure}
    \centering
    \begin{subfigure} {0.49\textwidth}
        \centering
        \includegraphics[width=\textwidth]{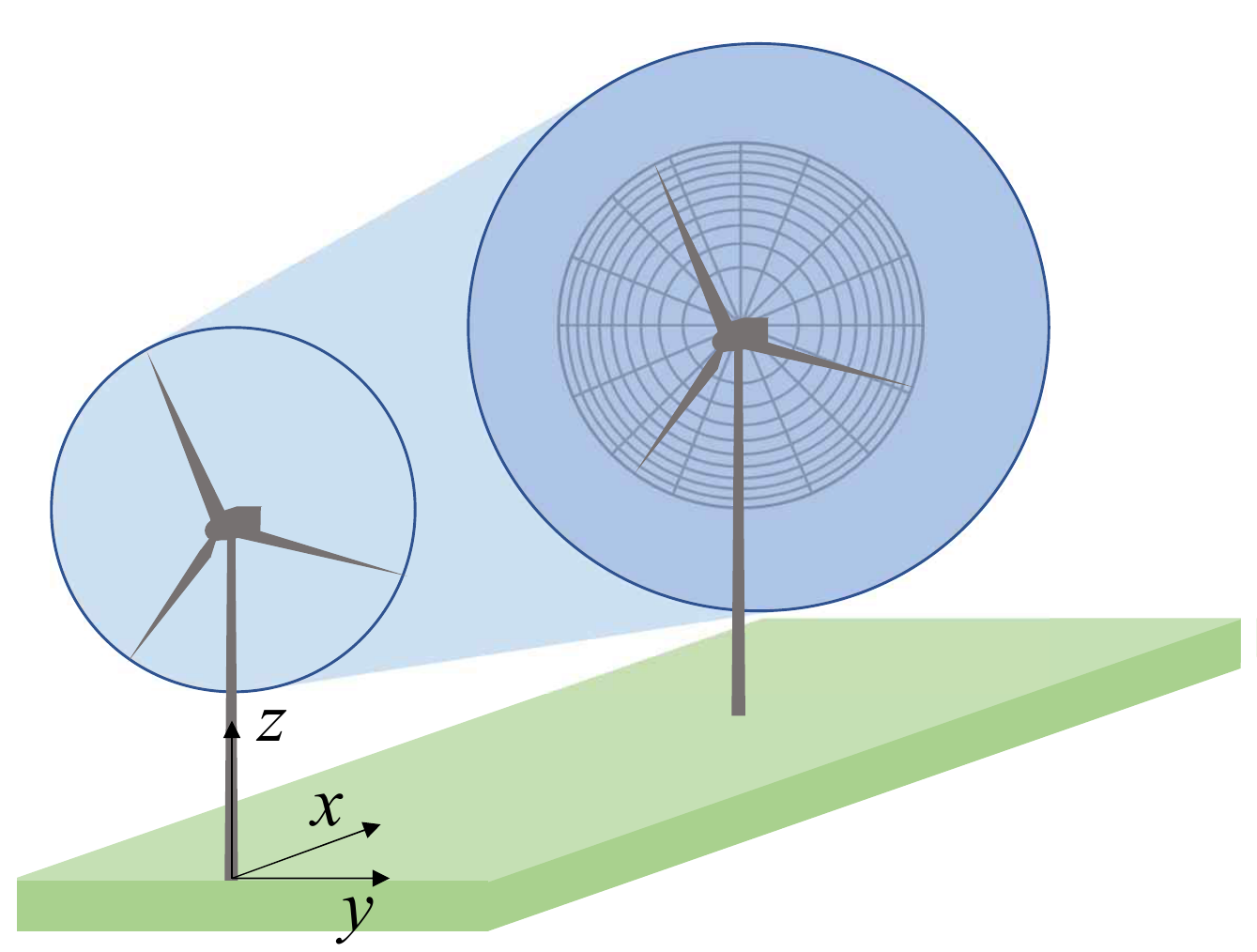}
        \caption{}
        \label{fig:unyawed}
    \end{subfigure}
    \hfill
    \begin{subfigure} {0.49\textwidth}
        \includegraphics[width=\textwidth]{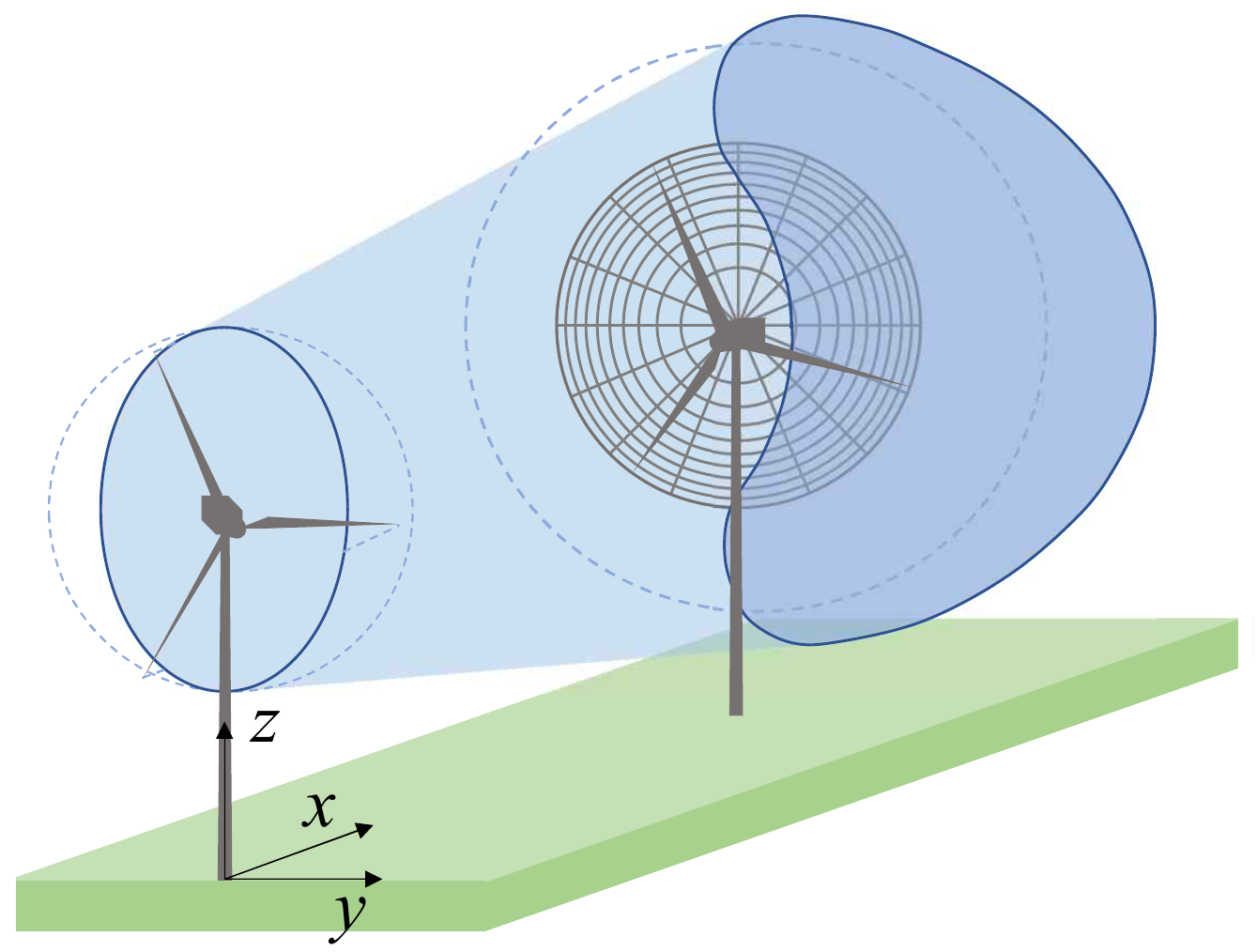}
        \caption{}
        \label{fig:yawed}
    \end{subfigure}	
    \caption{(a) Unyawed and (b) yawed wake shapes.}
    \label{fig:wake_model_yawed}
\end{figure}

    The offset of the wake center $Y_{\rm{offset}}$ increases proportionally with the square root of the streamwise distance $x$. Deformation of the wake center in the $y$ direction, $Y_{\rm{offset,z}}$ follows a Gaussian distribution in the $z$ direction. Fig.~\ref{fig:wake_model_yawed} shows the different wake shapes of the unyawed and the yawed wakes.

	\subsection{Wake merging model}\label{subsec:Wake merging model}
	
	Wakes by each turbine in the wind farm interfere with each other, therefore a wake merging model is needed. Four merging models are mainly used to represent wake merging phenomena as follows~\cite{lissaman1979, katic1986, voutsinas1990, niayifar2016}:
	\begin{gather}
	    u(x,y,z)=u_\infty(z)-\sum^n_{i=1}(u_\infty(z)-u_i(x,y,z)), \label{eq:lissaman1979} \\
	    u(x,y,z)=u_\infty(z)-\sqrt{\sum^n_{i=1}(u_\infty(z)-u_i(x,y,z))^2},\label{eq:katic1986} \\
	    u(x,y,z)=u_\infty(z)-\sqrt{\sum^n_{i=1}(\bar{u}_{in,i}-u_i(x,y,z))^2}, \label{eq:voutsinas1990} \\
		u(x,y,z)=u_\infty(z)-\sum^n_{i=1}(\bar{u}_{in,i}-u_i(x,y,z)), \label{eq:niayifar2016} 
	\end{gather}
	
    \noindent where $u(x,y,z)$ is the axial velocity field parallel to the direction of the atmospheric wind velocity, calculated with consideration of wake effects. $u_\infty(z)$ represents the wind speed in the atmosphere unaffected by the wind turbines. $u_i(x,y,z)$ means the wake speed of the $i$-th wind turbine. $\bar{u}_{in,i}$, which is the average wind speed entering the $i$-th wind turbine can be calculated as follows:
    \begin{equation} \label{windaverage}
	    \bar{u}_{in,i}=\frac{1}{A_i}\iint_{A_i} u(r,\theta) rdrd\theta.
	\end{equation}
    
    Eqs.~\eqref{eq:lissaman1979} and \eqref{eq:katic1986} use $u_\infty(z)$ to calculate the wake deficit caused by each wind turbine. This means that the distance between wind turbines is sufficiently far apart so that $\bar{u}_{in,i}$ is equal to $u_\infty(z)$. Eqs.~\eqref{eq:voutsinas1990} and \eqref{eq:niayifar2016} use $\bar{u}_{in,i}$ to calculate the wake deficit. This merging model is applicable even when the distance between wind turbines is relatively short. Eqs.~\eqref{eq:katic1986} and \eqref{eq:voutsinas1990} calculate the total wake deficit by summing the squares of the wake deficits of turbines. This method is called the energy deficit superposition. On the other hand, Eqs.~\eqref{eq:lissaman1979} and \eqref{eq:niayifar2016} calculate the total wake deficit by calculating the linear sum of the wake deficits of turbines. This method is called the velocity deficit superposition. It is shown that the velocity deficit superposition is more accurate than the energy deficit superposition from a large-eddy simulation study~\cite{niayifar2016}. In the present study, Eq.~\eqref{eq:niayifar2016} is used because it considers various distances between wind turbines based on the velocity deficit superposition.

\subsection{Optimization algorithm}\label{sec:optimization_algorithm}
 A single-objective constrained optimization problem can be formulated as
 \begin{equation} \label{eq:objectiveFunction}
    \min_{\bf{X}} f(\bf{X}), \quad \bf{X} \subset {\mathbb{R}}^{n \times m},
 \end{equation}
 subject to the constraints
  \begin{equation} \label{eq:constrainedFunction}
    g_i(\bf{X})\leq 0, \quad i=1,...,k,
 \end{equation}
 \noindent where $f(\bf{X})$ and $g(\bf{X})$ are cost and constraint functions for an $n \times m$ real variable matrix $\bf{X}$. When the number of dimensions exceeds 100, traditional evolutionary computation algorithms lose their effectiveness, which is usually known as the ``curse of dimensionality''~\cite{wang2021}. In this study, the algorithm for large-scale optimization is used because considering AYC in layout optimization greatly increases the dimension of $\bf{X}$. AGLDPSO can avoid suboptima with a high convergence speed~\cite{wang2020}.
 
 AGLDPSO is based on PSO~\cite{eberhart1995}. In PSO, each particle $j$ is updated at every step to find the optimal solution. Each particle has a position matrix $\bf{X}$ and a velocity matrix $\bf{V}$ which have the same dimension as $\bf{X}$. The searching strategy of PSO is described as follows:
 \begin{equation} \label{PSO_xv}
    \bf{X^{t+1}} = \bf{X^t}+\bf{V^t}, 
 \end{equation}
 where
 {

    $\bf{V^{t}} = \omega\times\bf{V^{t-1}} + c_1\times\bf{R_{1}^{t}}\odot(\bf{X_{pbest}^{t}}-\bf{X^{t}})+c_2\times\bf{R_{2}^{t}}\odot(\bf{X_{gbest}^{t}}-\bf{X^{t}}).$
    
 }
 \noindent Superscript t notates the present step, the symbol $\odot$ means element-wise product. $\bf{X_{pbest}}$ is the historical best position of $\bf{X}$ and $\bf{X_{gbest}}$ is the best position of $\bf{X_{pbest}}$. Velocities toward $\bf{X_{pbest}}$ and $\bf{X_{gbest}}$ are added to the current velocity of each particle. $\omega$ is the inertia considering the momentum of the previous movement. $c_1$ and $c_2$ are coefficients which can regulate the searching speed. $\bf{R_1}$ and $\bf{R_2}$ are random matrices.

Since all particles are updated toward $\bf{X_{gbest}}$ every iteration, it easily falls into a local optimum. To improve its performance, AGLDPSO uses a multi-subpopulation distributed model and an adaptive granularity learning strategy. The searching strategy of AGLDPSO is described as follows:
    \begin{equation} \label{AGLDPSO}
        \bf{X_w^{t+1}} = \bf{X_w^t}+\bf{V_w^t}, 
	\end{equation}
	where
	{
	
    	\centering{$\bf{V_w^t} = \bf{\Omega^t}\odot\bf{V_w^{t-1}}+c_1\times\bf{R_1^t}\odot(\bf{X_{sbest}^t}-\bf{X_w^t})+c_2\times\bf{R_2^t}\odot(\bf{X_{gbest}^t}-\bf{X_w^t}).$}
    	
	}

\begin{figure}
    \centering
    \begin{subfigure} {0.49\textwidth}
        \centering
        \includegraphics[width=\textwidth]{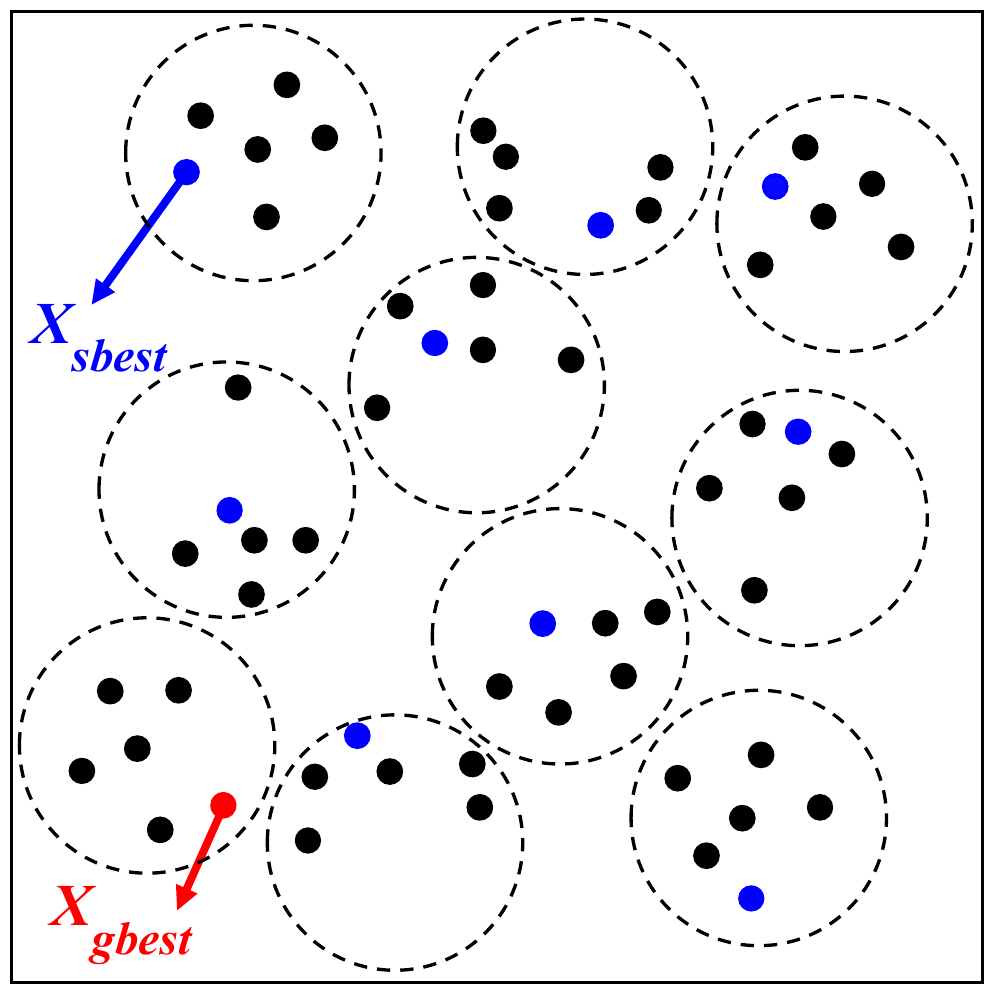}
        \caption{}
        \label{fig:MPSO_a}
    \end{subfigure}
    \hfill
    \begin{subfigure} {0.49\textwidth}
        \includegraphics[width=\textwidth]{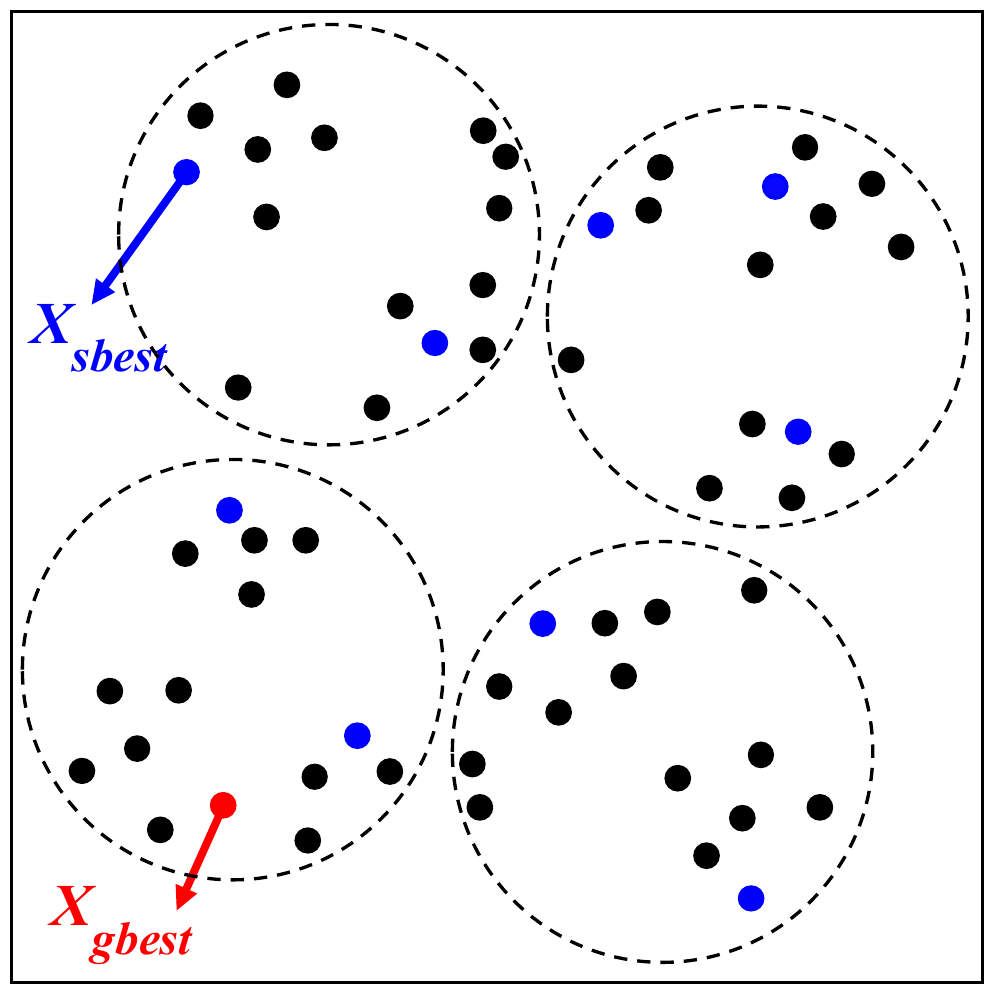}
        \caption{}
        \label{fig:MPSO_b}
    \end{subfigure}	
    \caption{Subpopulations of multi-swarm particle swarm optimization for different subpopulation sizes. (a) Small size and (b) large size.}
    \label{fig:MPSO}
\end{figure}

\noindent Subscript $w$ notates the worst particle in each subpopulation, and $\bf{\Omega}$ is random matrix. $\bf{X_{sbest}}$ is the best position of each subpopulation and $\bf{X_{gbest}}$ is the best position of $\bf{X_{sbest}}$. In Eq.~\eqref{AGLDPSO}, only the worst particle in each subpopulation is updated to search for the global optimal point. More subpopulations lead to more exploration. In each iteration, AGLDPSO divides the particles into subpopulations of the automatically calculated size. For example, when most particles are close to the global worst particle, the size of the subpopulation becomes smaller to increase the number of subpopulations, and the exploration rate increases. On the other hand, fewer subpopulations lead to less exploration. In Fig.~\ref{fig:MPSO}, because the number of subpopulations in Fig.~\ref{fig:MPSO}(\subref{fig:MPSO_a}) is more than that in Fig.~\ref{fig:MPSO}(\subref{fig:MPSO_b}), the situation in Fig.~\ref{fig:MPSO}(\subref{fig:MPSO_a}) has a tendency to explore rather than exploit. These mechanisms enhance convergence to an optimal point.
  
\section{Problem setup}\label{sec:Problem setup}
    \subsection{AEP calculation}\label{sec:AEP_calculation}
    
\begin{figure}
    \centering
    \begin{subfigure} {0.49\textwidth}
        \centering
        \includegraphics[width=\textwidth]{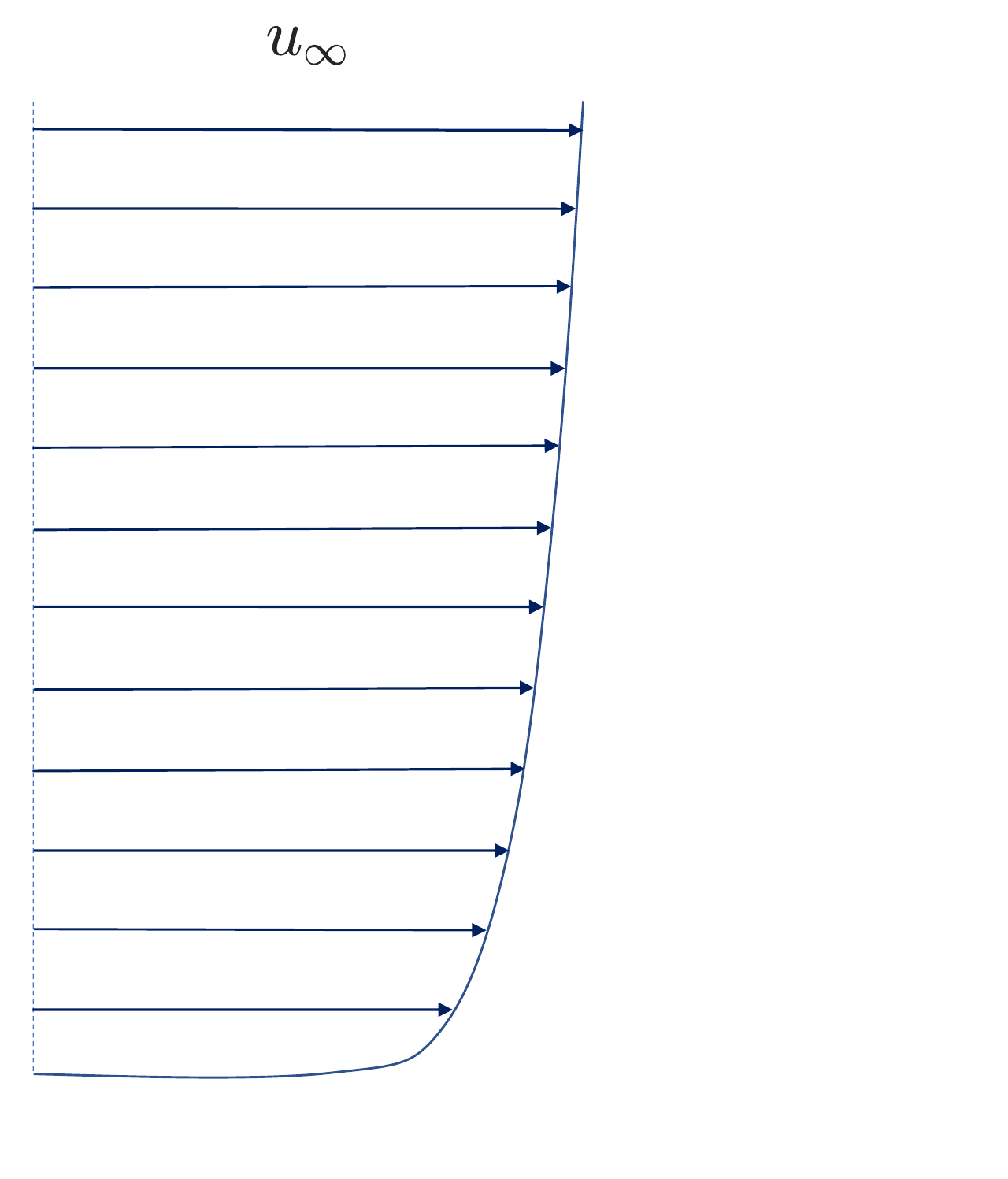}
        \caption{}
        \label{fig:ABL}
    \end{subfigure}
    \hfill
    \begin{subfigure} {0.49\textwidth}
        \includegraphics[width=\textwidth]{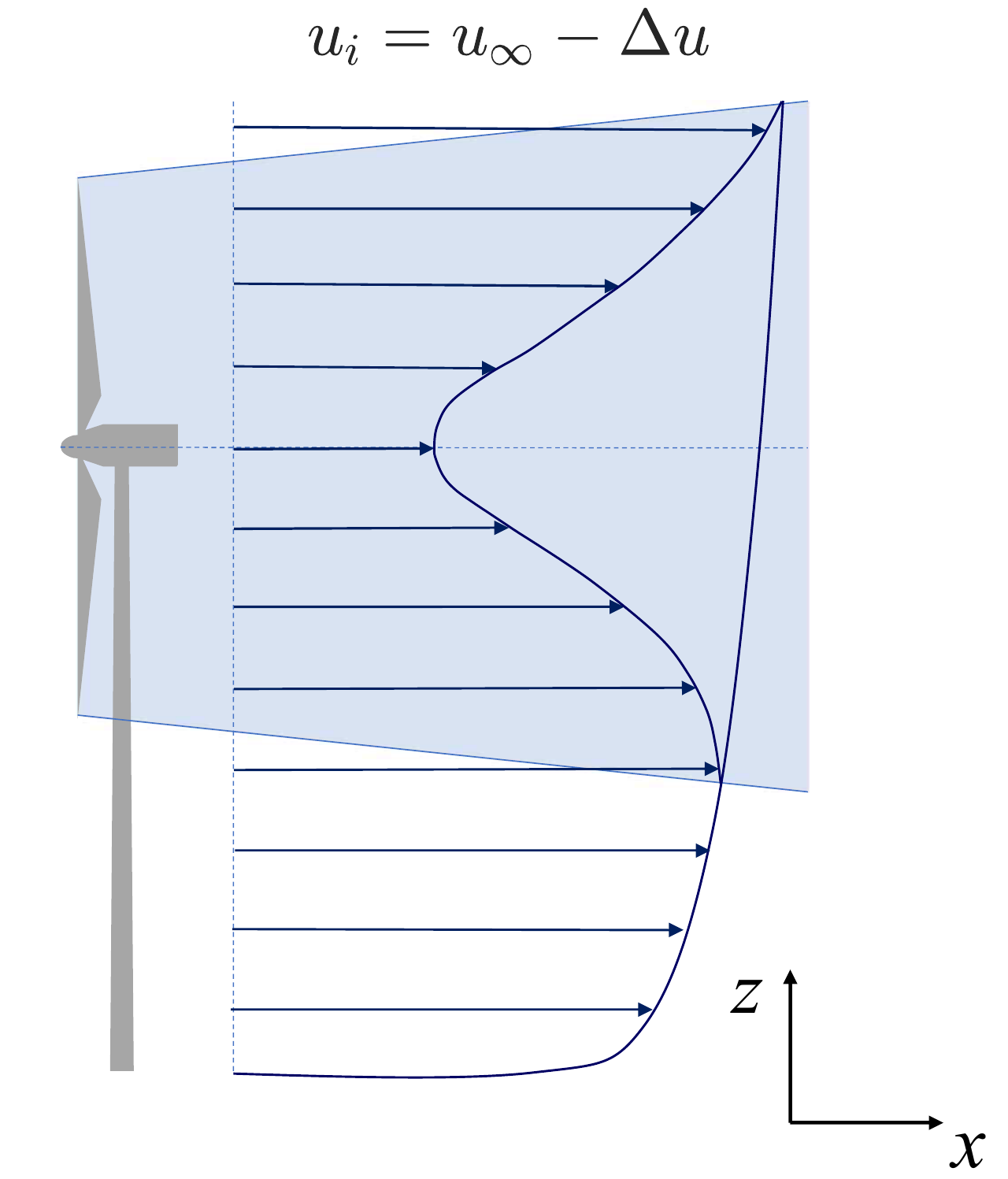}
        \caption{}
        \label{fig:ABL_with_wake}
    \end{subfigure}	
    \caption{(a) The atmospheric boundary layer profile and (b) the wake velocity profile.}
    \label{fig:ABL_wake}
\end{figure}

    The free stream velocity $u_\infty(z)$ is calculated using an atmospheric boundary layer equation. The power law equation of an atmospheric boundary layer is as follows:
    \begin{equation} \label{eq:ABL}
	    u_\infty(z) = U_{ref}(z/z_{ref})^\alpha,
	\end{equation}
	
	\noindent where $U_{ref}$ is the measured wind speed at the reference height $z_{ref}$ and $\alpha$ is an exponent determined empirically. In this study, $U_{ref}=$ 6m/s or 8m/s, and $z_{ref}=$ 25m. $\alpha$ is set to 0.1 for an offshore wind farm~\cite{wang2018}. Fig.~\ref{fig:ABL_wake} shows the profile of the atmospheric boundary layer and the wind turbine wake.
    
\begin{figure}
    \centering
    \includegraphics[width=\textwidth]{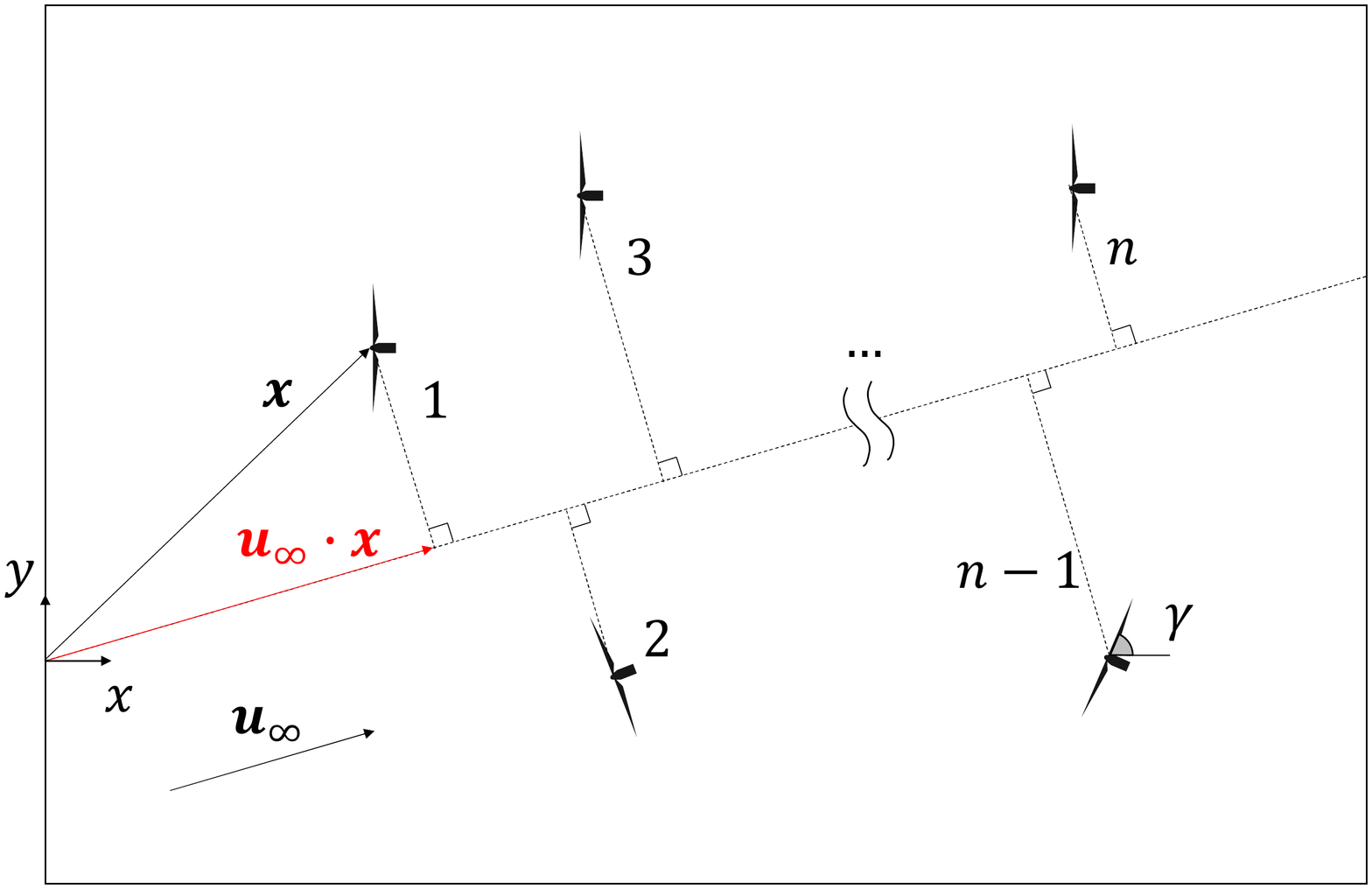}
    \caption{Schematic illustration of the method sorting wind turbines along the wind direction.}
    \label{fig:sorting}
\end{figure}

	Using $u_\infty(z)$ and the wake merging model, $\bar{u}_{in,i}$ of Eq.~\eqref{windaverage} can be obtained. Because of wake interactions, $\bar{u}_{in,i}$ must be calculated on wind turbines along the wind direction. Fig.~\ref{fig:sorting} explains how to sort wind turbines along the wind direction. The inner product of the wind vector and the position vector of a wind turbine $\bf{u_\infty}\cdot\bf{x}$ is calculated to sort wind turbines. Numerical integration in Eq.~\eqref{windaverage} is conducted using the grid of the wind turbine as represented in Fig.~\ref{fig:wake_model_yawed}.
	
	The power coefficient $C_P$ and the thrust coefficient $C_T$ for each wind turbine are functions of $\bar{u}_{in,i}$. They depend on the type of the wind turbine. Fig.~\ref{fig:Vestas} shows the power curves, the power coefficients, and the thrust coefficients of the Vestas V80 and the Vestas V112 models~\cite{v80_wu2012, v112_amin2021}. The generated power of the $i$-th wind turbine is calculated as follows: 
	\begin{equation} \label{power}
	    P = \sum_i \frac{1}{2}\rho C_P(\bar{u}_{in,i}) \bar{u}_{in,i}^3 A_i.
	\end{equation}
	
	In order to calculate AEP of a wind farm, a probability distribution function $p(U_{\infty,j})$ of the annual wind vector is required. Then, AEP can be obtained as follows:
	\begin{equation} \label{aep}
	    AEP = 365\times24\times\sum_j p(U_{\infty,j})P|_{u_\infty=U_{\infty,j}},
	\end{equation}
	
	\noindent where $U_{\infty,j}$ is the $j$-th wind vector of $p(U_{\infty,j})$. Units of $P$ and $AEP$ are [W] and [Wh], respectively.
	
    \subsection{Optimization configurations}\label{sec:optimization_configurations}
	To solve the optimization problem, an appropriate cost function $C(\bf{X})$, constraints, and a variable matrix $\bf{X}$ have to be defined. For the joint optimization, which consists of AYC and layout optimization, the cost function is set as the reciprocal of AEP of a wind farm represented as follows:
	\begin{equation} \label{joint_def}
        [\bf{X_{pos}^{joint}}, \bf{X_{yaw}^{joint}}]=\underset{\bf{X}}{\text{argmin}}\ C(\bf{X}), \quad \bf{X} \subset {\mathbb{R}}^{n \times (m+2)},
	\end{equation}
	where
	{
	    \[\bf{X} = [\bf{X_{pos}}, \bf{X_{yaw}}],\]
	    
	    \[\bf{X_{pos}}=
	    \left[\begin{matrix}
            x_1 & y_1\\
            x_2 & y_2\\
            \vdots & \vdots\\
            x_n & y_n\\
        \end{matrix}\right],\quad \bf{X_{yaw}}=
	    \left[\begin{matrix}
            \gamma_{11} & \gamma_{12} & \dots & \gamma_{1m}\\
            \gamma_{21} & \gamma_{22} & \dots & \gamma_{2m}\\
            \vdots & \vdots & \ddots & \vdots\\
            \gamma_{n1} & \gamma_{n2} & \dots & \gamma_{nm}\\
        \end{matrix}\right],\]
        
        \[C(\bf{X}) = \frac{1}{AEP|_{\bf{X}}},\]
        
	}
	
	\noindent subject to constraints
	\begin{gather}
	    x_{min}\leq{x_i}\leq{x_{max}}, \quad \forall i, \label{eq:x} \\
	    y_{min}\leq{y_i}\leq{y_{max}}, \quad \forall i, \label{eq:y} \\
	    (x_i-x_j)^2+(y_i-y_j)^2>D^2, \quad \forall i,j. \label{eq:distance}
	\end{gather}
	
	\noindent The number of turbines in a wind farm is $n$, and the wind rose consists of $m$ winds. $\bf{X}$ consists of positions and yaw angles of $n$ wind turbines. $AEP|_{\bf{X}}$ is AEP of a wind farm where the variable matrix is given as $\bf{X}$. Constraints are defined in Eqs.~(\ref{eq:x})--(\ref{eq:distance}). Wind turbine positions have to be optimized in the wind farm domain [$x_{min}$, $x_{max}$]$\times$[$y_{min}$, $y_{max}$], and the distances between each wind turbine have to be more than the wind turbine diameter $D$. The present joint optimization algorithm finds the optimal layout [$\bf{X_{pos}^{joint}}$, $\bf{X_{yaw}^{joint}}$] where AEP is the maximum using AGLDPSO.
	
	To analyze the joint optimization performance against the conventional layout optimization, the sequential optimization is defined as follows:
	\begin{equation} \label{sequent_greedy}
        \bf{X_{pos}^{seq}}=\underset{\bf{X_{pos}}}{\text{argmin}}\ C([\bf{X_{pos}},\bf{\Theta}]),
	\end{equation}
		\begin{equation} \label{sequent_AYC}
        \bf{X_{yaw}^{seq}}=\underset{\bf{X_{yaw}}}{\text{argmin}}\ C([\bf{X_{pos}^{seq}},\bf{X_{yaw}}]).
	\end{equation}
	
	\noindent The conventional layout optimization for finding an optimal layout using a greedy yaw control strategy $\bf{\Theta}$ is first performed. Using its optimal layout $\bf{X_{pos}^{seq}}$, the yaw angles of wind turbines according to $m$ wind directions are optimized to maximize AEP.

	\section{Results and discussions}\label{sec:Optimization result}
	
	\subsection{Comparison of optimization methods }\label{subsec:Optimization methods comparison}
	
	The result of the present method is compared to that of DBHM using the WF1 case which is shown in Fig.~\ref{fig:WFLO_cases}(\subref{fig:5by5_5D_uniform}). The domain of the WF1 case is defined as [0m, 1600m]$\times$[0m, 1600m] and the layout of 25 Vestas V80s is optimized. Wind which has the same speed and the same probability blows from eight directions. $U_{ref}$ is 8 m/s, so $\bar{u}_{in,1}$ in Eq.~\eqref{windaverage} is 8.85 m/s for the Vestas V80. For comparison of results,  $AEP|_{\bf{X}=[\bf{X_{pos}^{seq}}, \Theta]}$,
	$AEP|_{\bf{X}=[\bf{X_{pos}^{seq}}, \bf{X_{yaw}^{seq}}]}$,
	$AEP|_{\bf{X}=[\bf{X_{pos}^{joint}}, \Theta]}$, and
	$AEP|_{\bf{X}=[\bf{X_{pos}^{joint}}, \bf{X_{yaw}^{joint}}]}$ are defined as $S_\Theta$, $S_{AYC}$, $J_{\Theta}$, and $J_{AYC}$, respectively. $S_{\Theta}$ indicates AEP of conventional layout optimization based on a greedy control strategy, and $J_{AYC}$ is AEP of the joint optimization.
	
\begin{figure}
	\centering
	\begin{subfigure} {0.32\textwidth}
        \includegraphics[width=\textwidth]{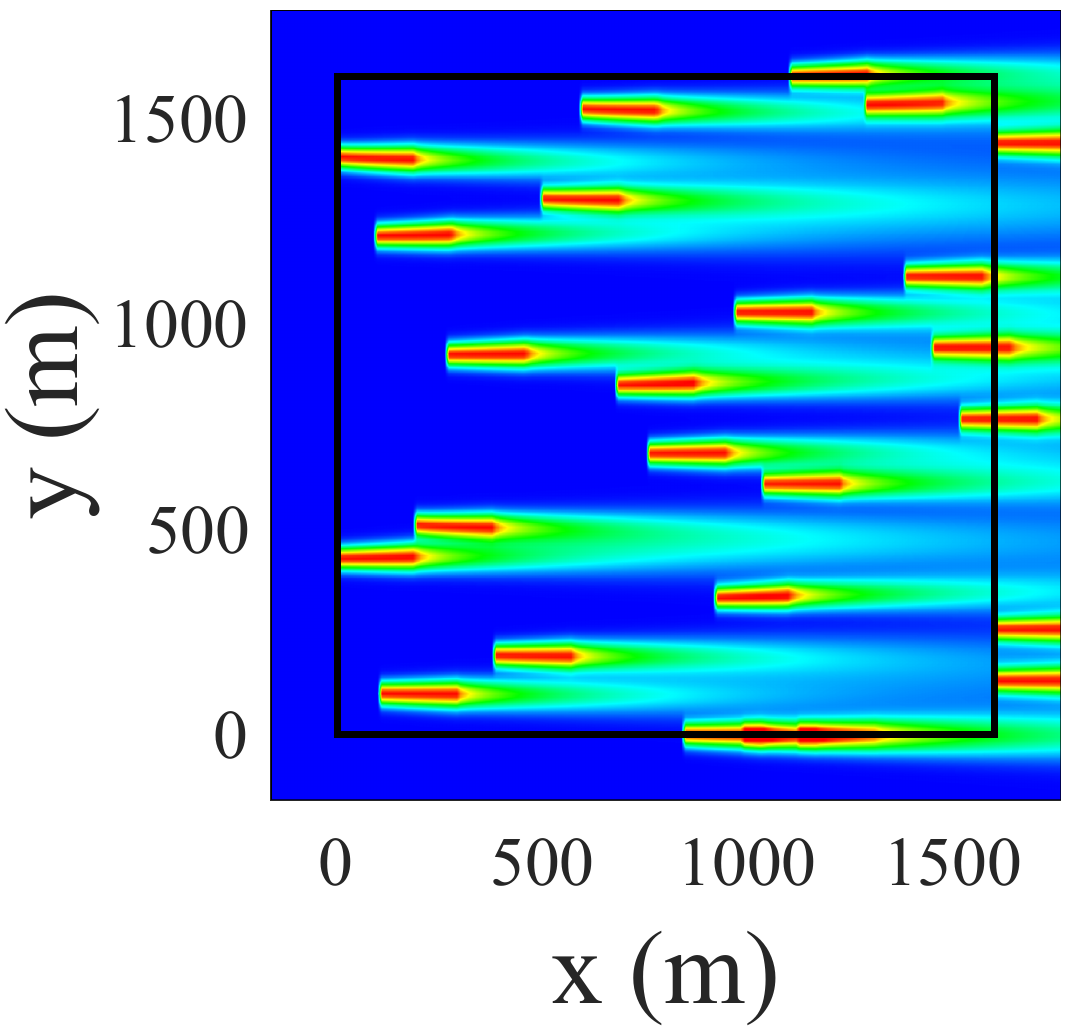}
        \caption{}
        \label{fig:5by5_S_AYC1}
    \end{subfigure}
    \hfill
    \begin{subfigure} {0.32\textwidth}
        \includegraphics[width=\textwidth]{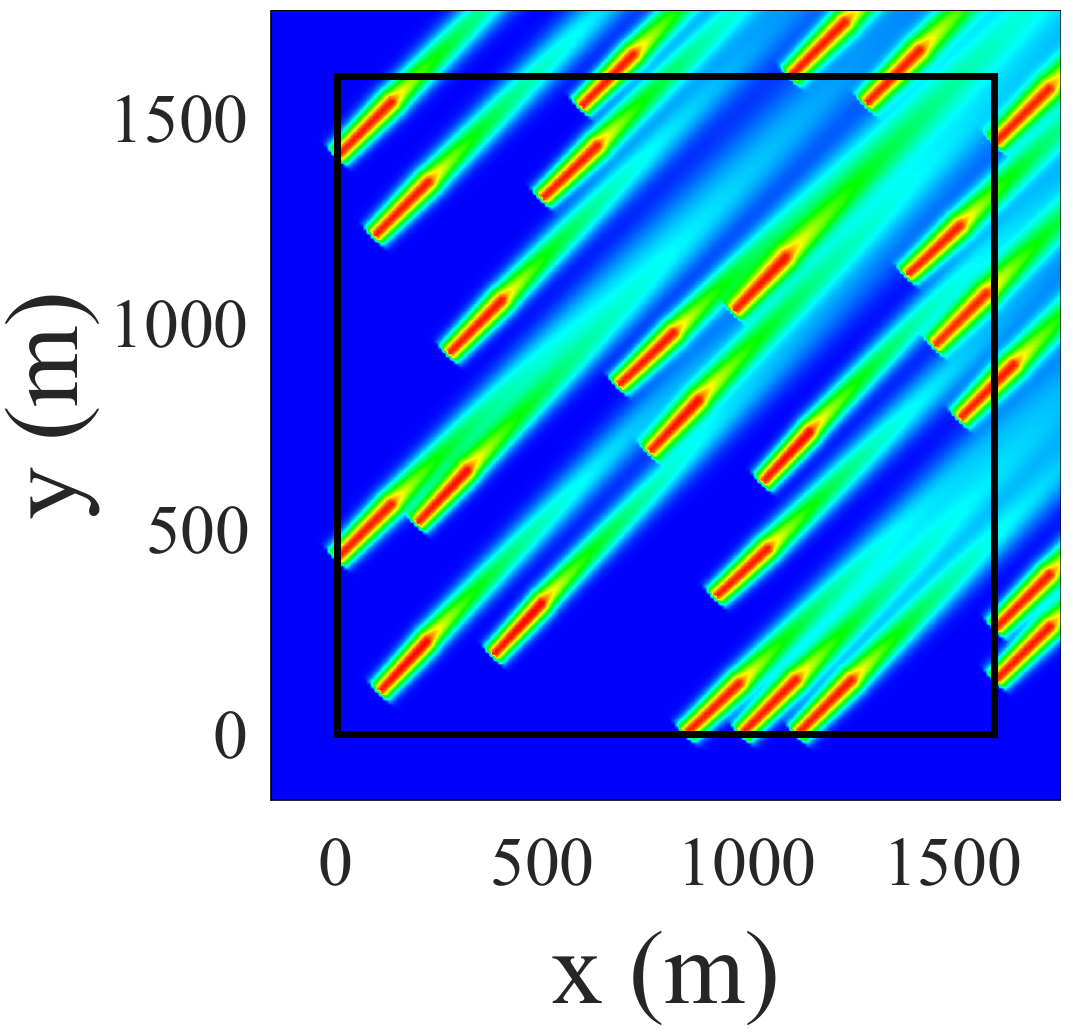}
        \caption{}
        \label{fig:5by5_S_AYC2}
    \end{subfigure}
    \hfill
    \begin{subfigure} {0.32\textwidth}
        \includegraphics[width=\textwidth]{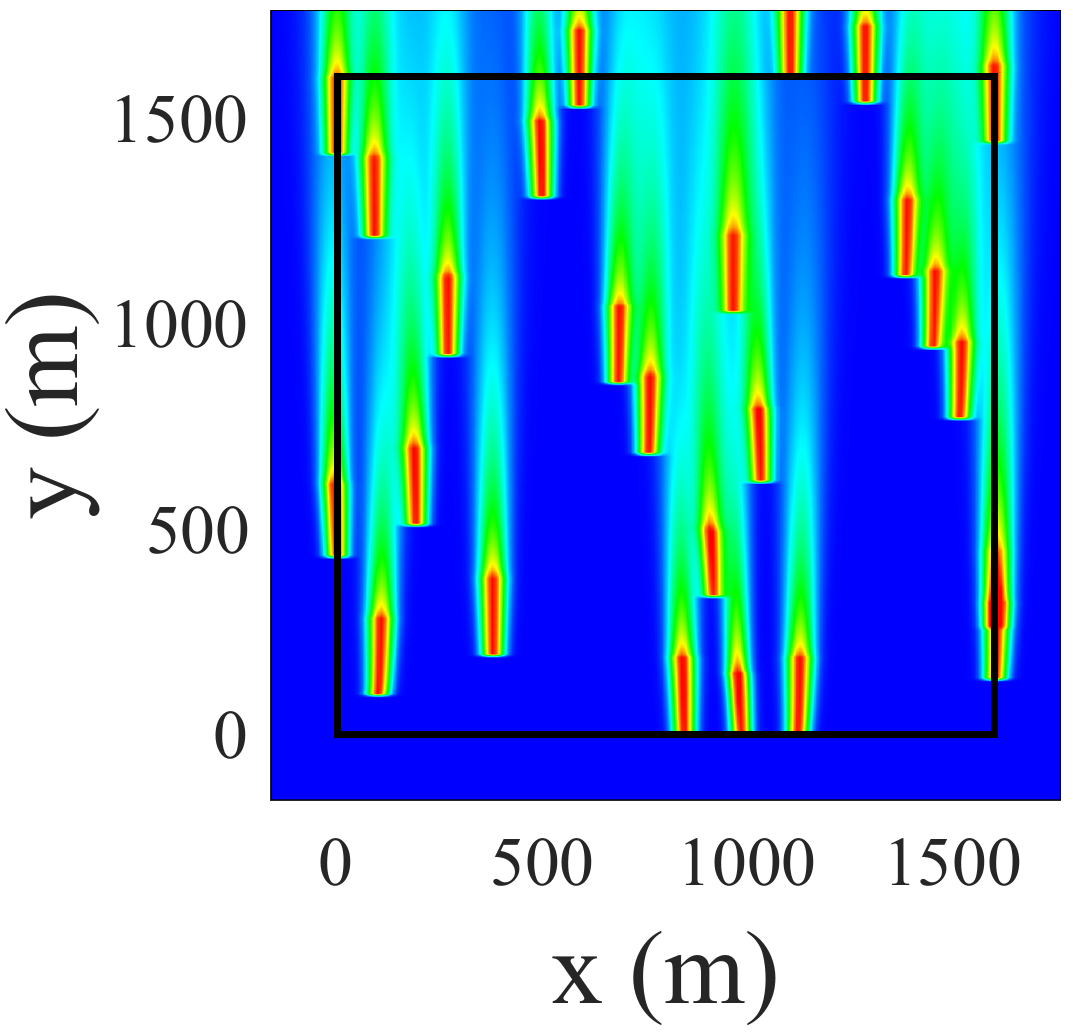}
        \caption{}
        \label{fig:5by5_S_AYC3}
    \end{subfigure}

    \begin{subfigure} {0.32\textwidth}
        \includegraphics[width=\textwidth]{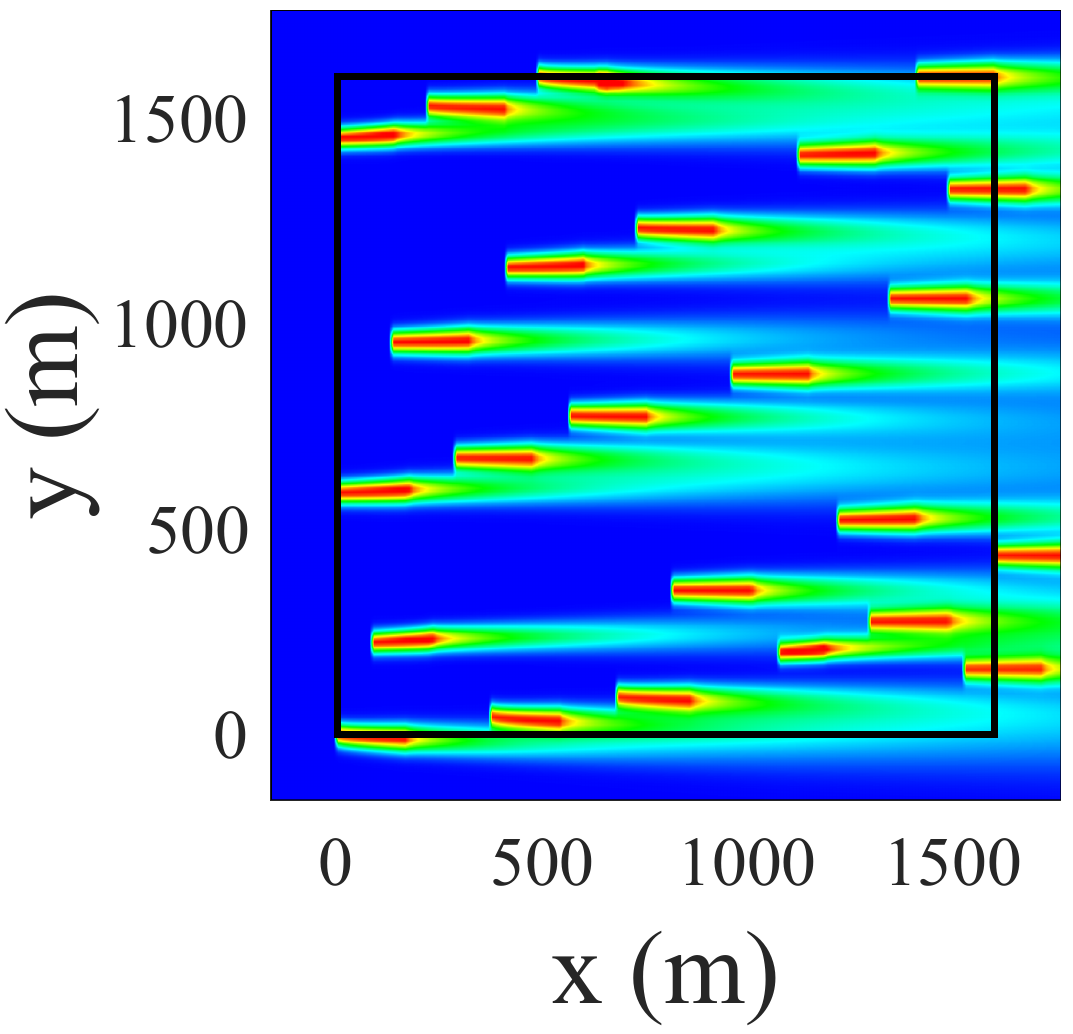}
        \caption{}
        \label{fig:5by5_J_AYC1}
    \end{subfigure}
    \hfill
    \begin{subfigure} {0.32\textwidth}
        \includegraphics[width=\textwidth]{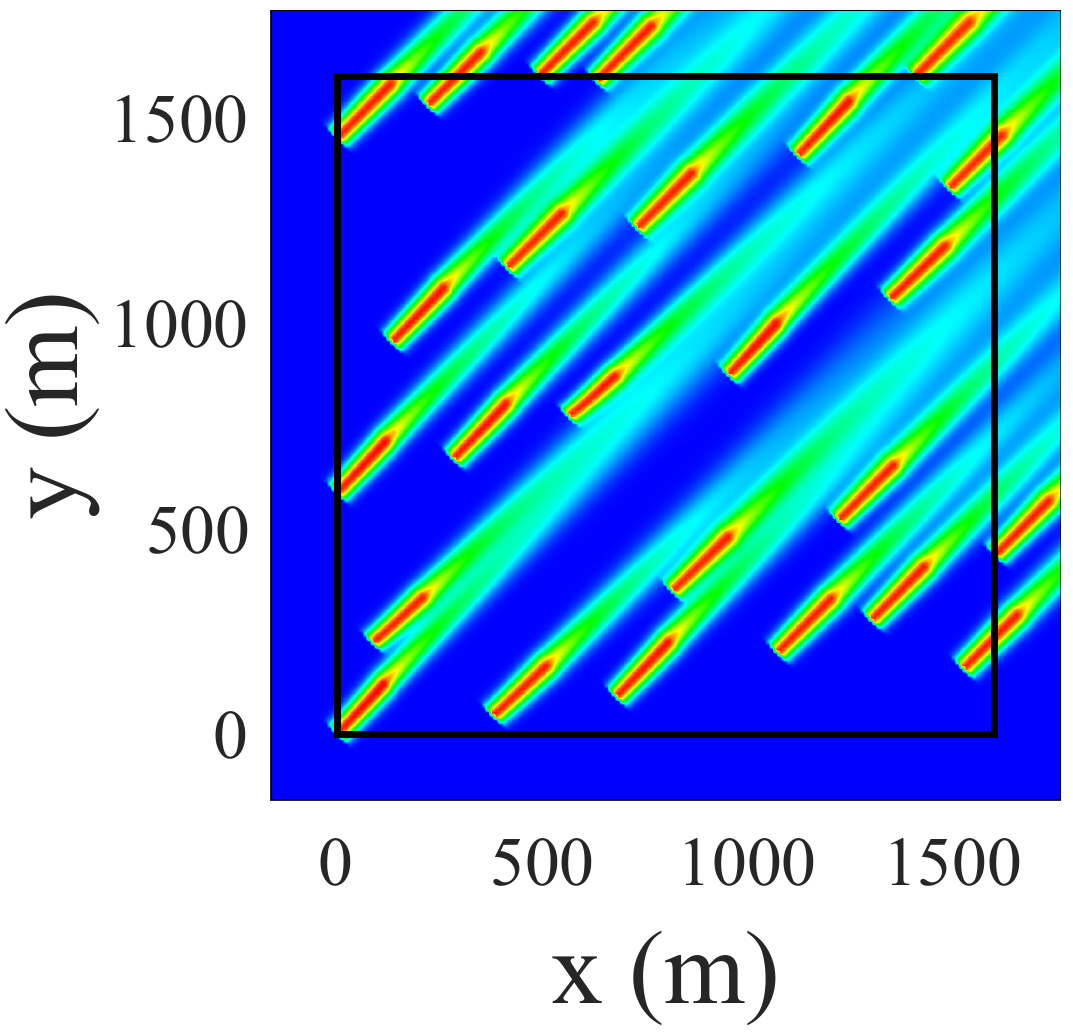}
        \caption{}
        \label{fig:5by5_J_AYC2}
    \end{subfigure}
    \hfill
    \begin{subfigure} {0.32\textwidth}
        \includegraphics[width=\textwidth]{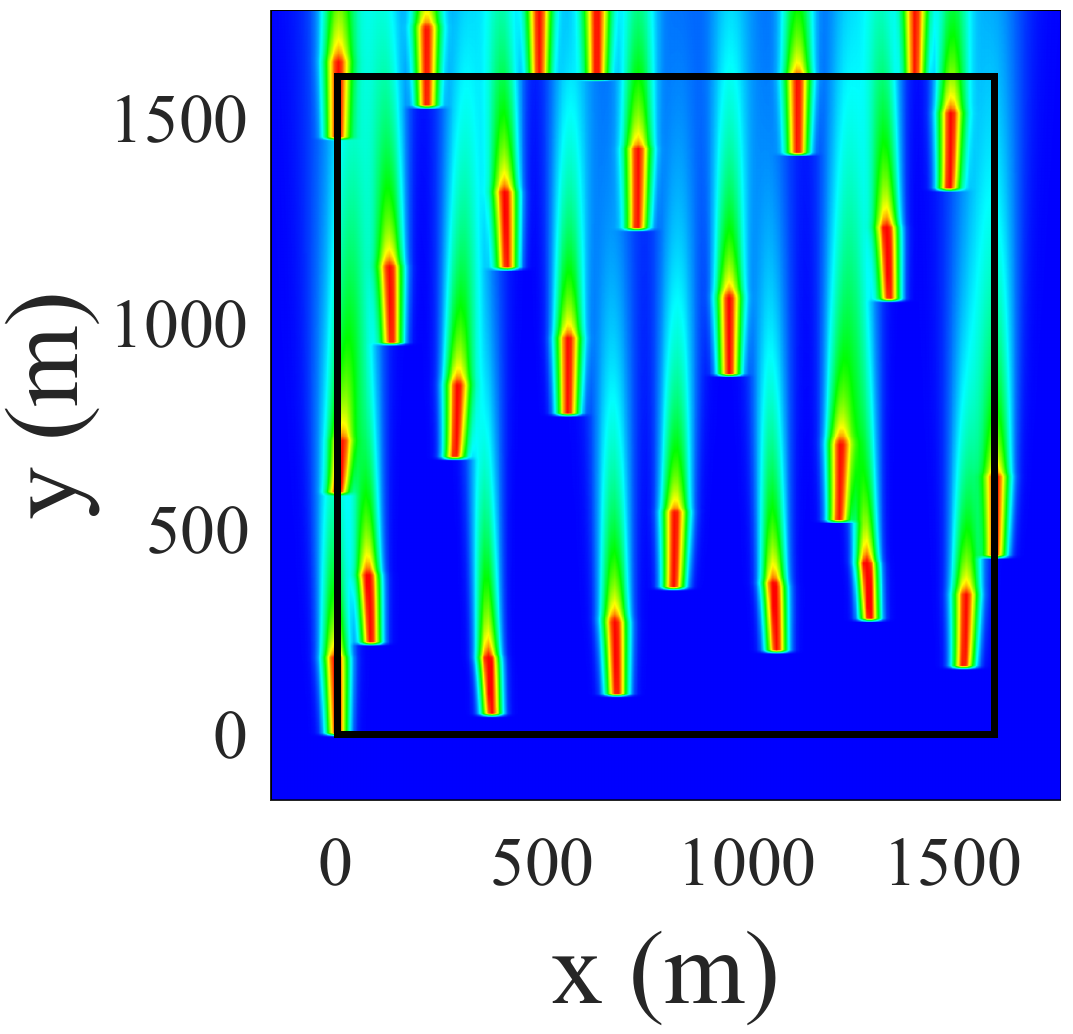}
        \caption{}
        \label{fig:5by5_J_AYC3}
    \end{subfigure}
    \begin{subfigure} {0.6\textwidth}
        \includegraphics[width=\textwidth]{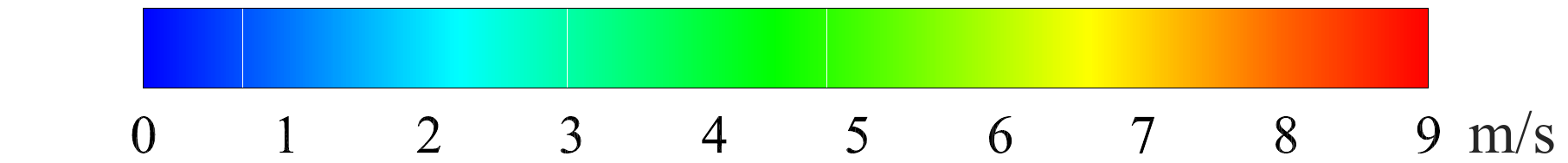}
        \label{fig:colorbar_ms}
    \end{subfigure}
	\caption{Contours of wakes of sequentially optimized layouts for the wind from (a) the west, (b) the southwest, and (c) the south. Contours of wakes of jointly optimized layouts for the wind from (a) the west, (b) the southwest, and (c) the south. The color bar corresponds the magnitude of the wake deficit.}
	\label{fig:combined_5by5_5D}
\end{figure}

	Optimized layouts of the WF1 case using the sequential optimization method and the present method are shown in Fig.~\ref{fig:combined_5by5_5D}. In both layouts, wakes of the upstream wind turbines are deflected by AYC. Downstream wind turbines are located as far from the wake region as possible. The difference in power production between greedy control and AYC is shown in Table~\ref{tab:5by5_5D}. Improvements of $S_{AYC}$, $J_{\Theta}^{DBHM}$, $J_{AYC}^{DBHM}$, $J_{\Theta}^{AGLD}$, and $J_{AYC}^{AGLD}$ compared to $S_{\Theta}$ are shown. Superscripts $DBHM$ and $AGLD$ represent that the results are optimized using the DBHM and the AGLDPSO algorithms, respectively. $S_{\Theta}$ is AEP of the optimized layout using greedy control $\Theta$, so $J_\Theta$ has to be lower than $S_{\Theta}$. It means that the optimal layout of the joint optimization problem is the solution by AYC, not by greedy control. AEP of AYC is always higher than AEP of greedy control and these improvements are higher for the jointly optimized layout than the sequentially optimized layout. As a result, $J_{AYC}^{AGLD}$ is 3.73\% higher than $S_{\Theta}$. Compared to $S_{AYC}$, $J_{AYC}^{AGLD}$ is 2.42\% higher. To maximize the efficiency of AYC, the joint optimization problem must be solved at the design stage of a wind farm layout.
	
	In Table~\ref{tab:5by5_5D}, $J_{AYC}^{DBHM}$ is 2.08\% higher than $S_\Theta$, but it is 1.62\% lower than $J_{AYC}^{AGLD}$. The main differences between the DBHM method and the present method are the initialization of $\bf{X}$ and the searching strategy. The DBHM method starts with an optimal layout of the sequential optimization problem and uses the optimal layouts of decomposed subproblems at every optimization step. The solution of the DBHM method depends on the initial layout $\bf{X_{pos}^{seq}}$. On the other hand, The present method searches the globally optimal solution for all wind velocities among particles which start from numerous randomly distributed initial layouts. As a result, $||\bf{X_{pos}^{joint}}-\bf{X_{pos}^{seq}}||_2$ of the DBHM method is found to be 59 m, while that of the present method is 1580 m.
	
	The optimal layout of the joint optimization problem can be very different from that of the sequential optimization problem. The present method searches for a globally optimal layout of joint optimization and does not depends on an initial layout. Therefore, the present method is more accurate for large-scale joint optimization than the DBHM method, which starts from an initial layout. In this study, the remaining case results are obtained using the present method.
	
	\subsection{Power curve gradients of wind turbines}\label{subsec:Power curve gradient}
	
	Commercial wind turbines have their own characteristics, such as power efficiency and design parameters. These differences affect the optimal layout of the joint optimization problem and the improvement of $J_{AYC}$ compared to $S_{AYC}$. To investigate the effects of the wind speed and the wind turbine model on joint optimization performance, two modified WF1 cases are considered. For the first case, $U_{ref}$ is changed from 8 m/s to 6 m/s, and for the second case, the wind turbine model is changed from the Vestas V80 to the Vestas V112. The results are shown in Table~\ref{tab:modifiedWF1}.
	
\begin{figure}
    \centering
    \begin{subfigure} {0.32\textwidth}
        \includegraphics[width=\textwidth]{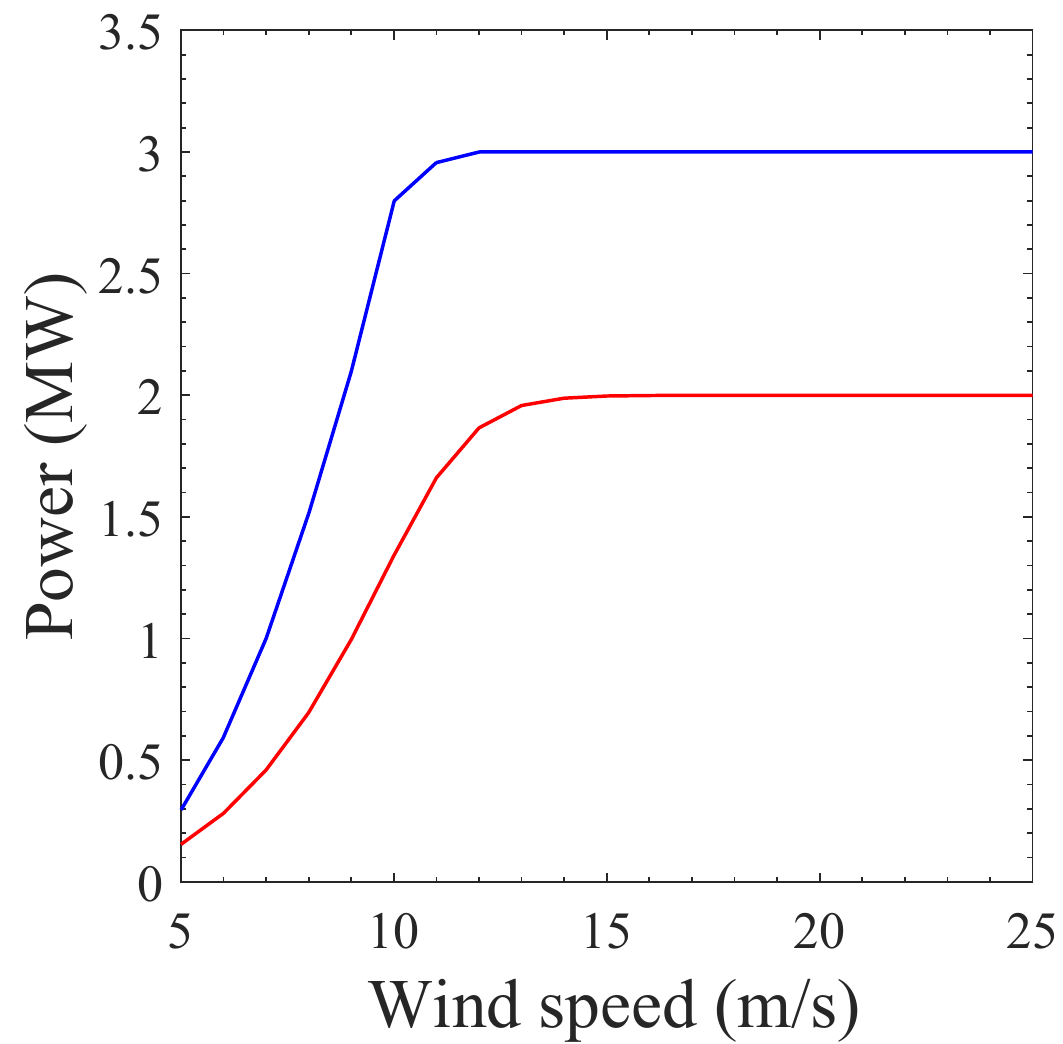}
        \caption{}
        \label{fig:power}
    \end{subfigure}
    \hfill
    \begin{subfigure} {0.32\textwidth}
        \includegraphics[width=\textwidth]{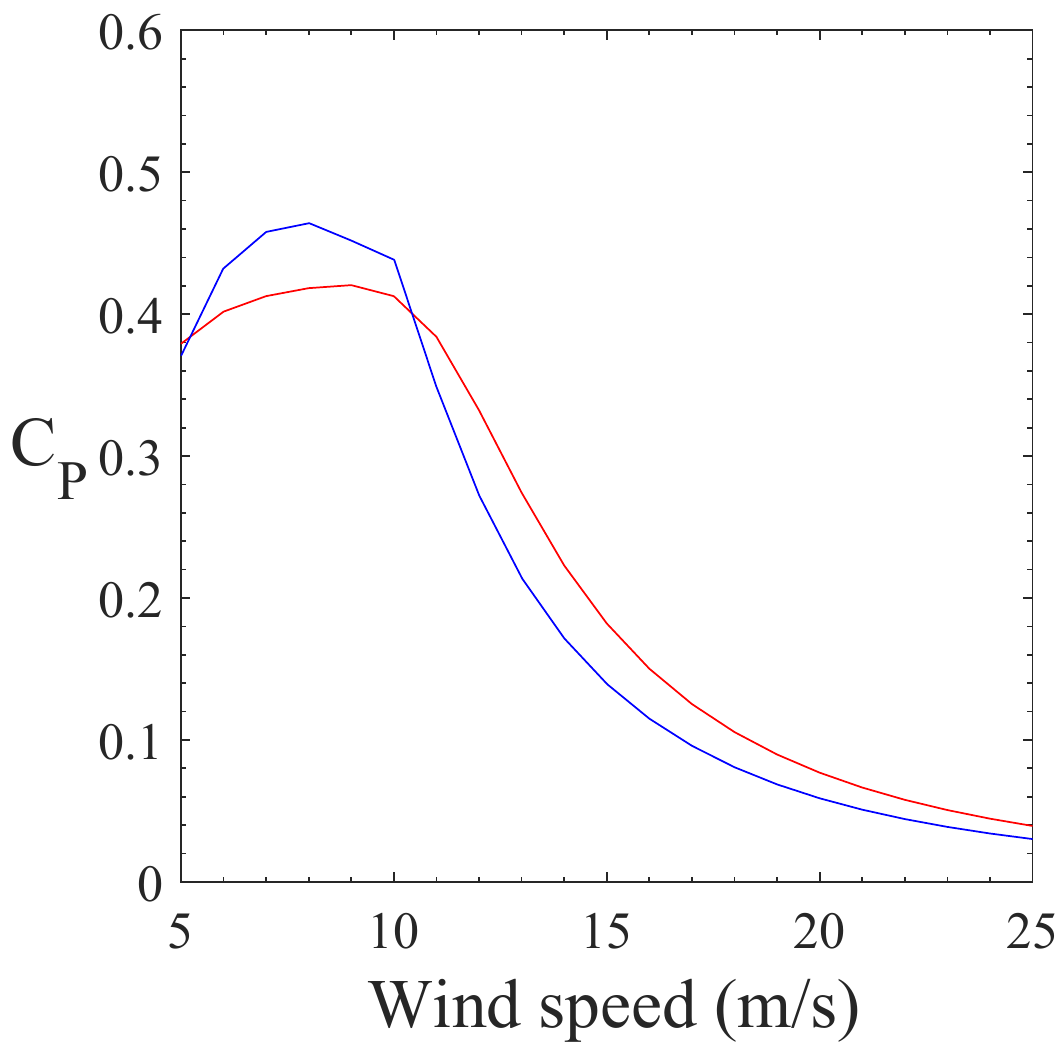}
        \caption{}
        \label{fig:C_P}
    \end{subfigure}
    \hfill
    \begin{subfigure} {0.32\textwidth}
        \includegraphics[width=\textwidth]{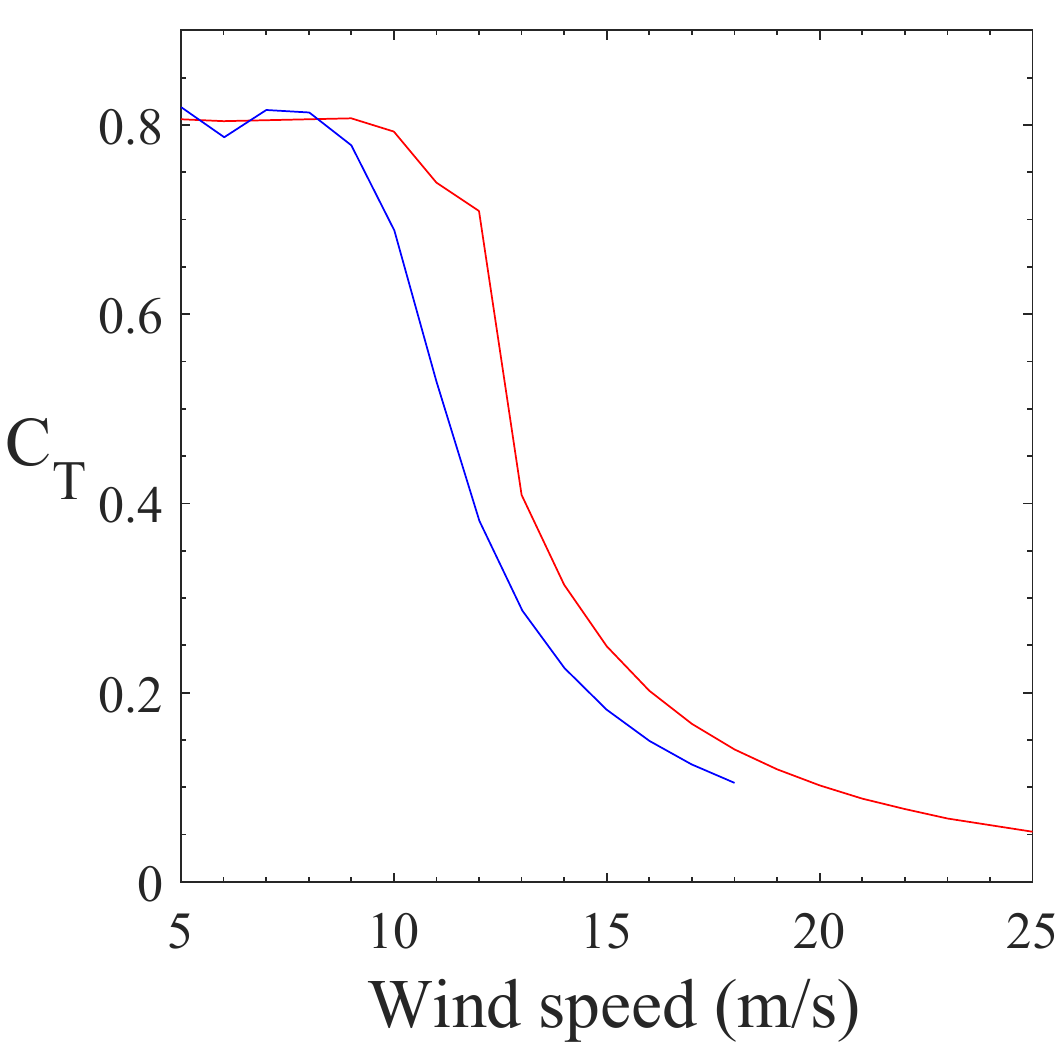}
        \caption{}
        \label{fig:C_T}
    \end{subfigure}	
    \caption{(a) Power curves, (b) power coefficients, and (c) thrust coefficients of the Vestas V80 (\red{\sampleline{}}) and the Vestas V112 (\blue{\sampleline{}}), respectively.}
    \label{fig:Vestas}
\end{figure}

\begin{figure}
    \centering
    \begin{subfigure} {0.49\textwidth}
        \includegraphics[width=\textwidth]{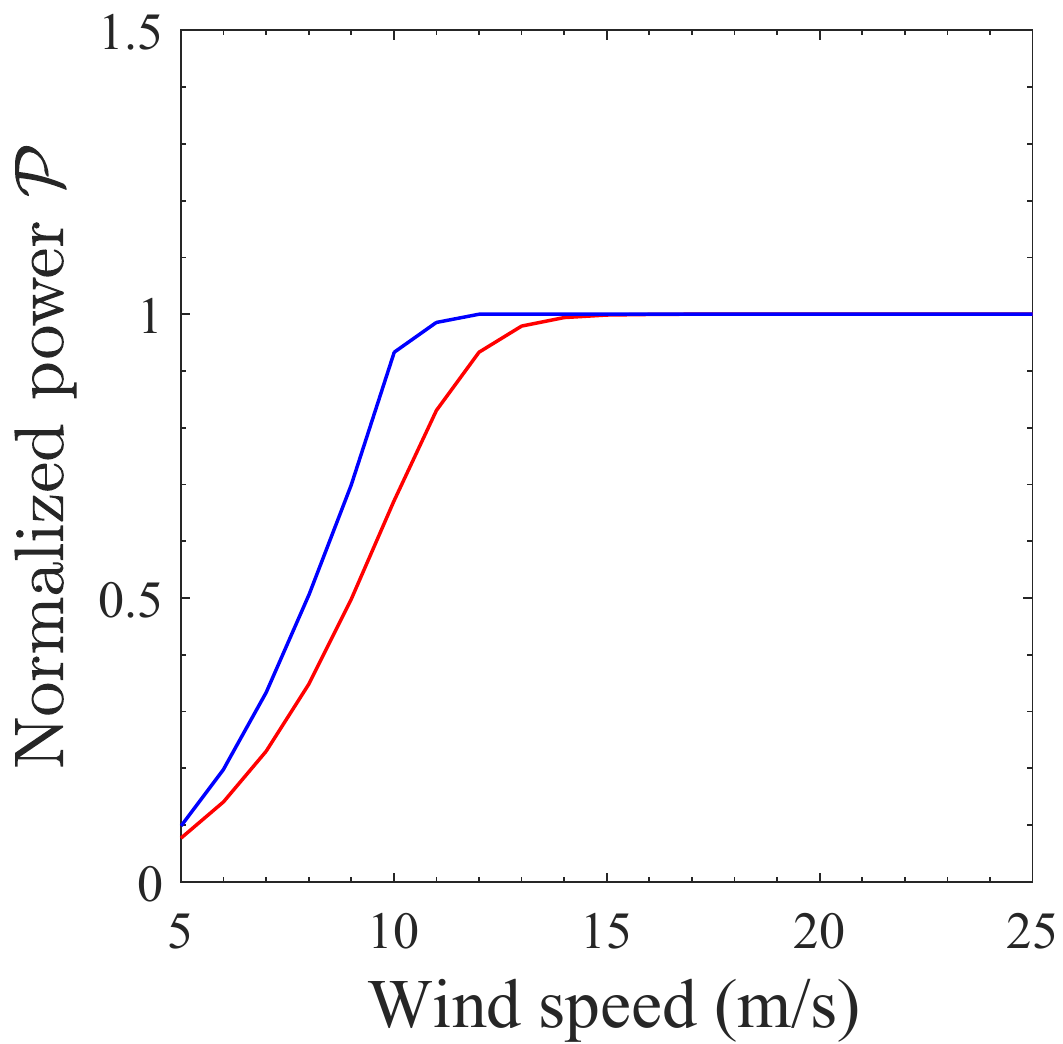}
        \caption{}
        \label{fig:normalpower}
    \end{subfigure}
    \hfill
    \begin{subfigure} {0.49\textwidth}
        \includegraphics[width=\textwidth]{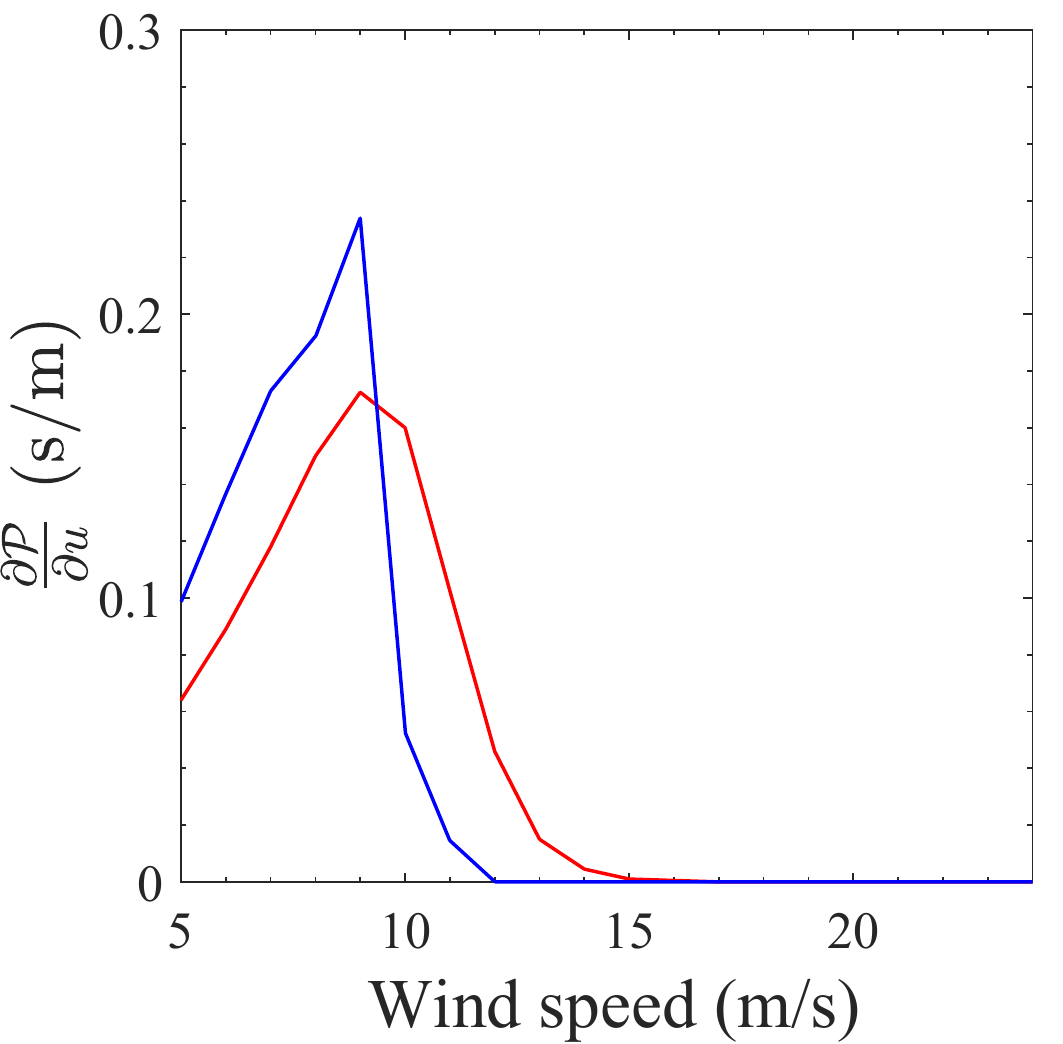}
        \caption{}
        \label{fig:power_derivative}
    \end{subfigure}
	\caption{(a) Power curves of wind turbines normalized by the rated power and (b) derivatives of the power curve with respect to wind speed. Wind turbines are the Vestas V80 (\red{\sampleline{}}) and the Vestas V112 (\blue{\sampleline{}}).}
	\label{fig:derivative}
\end{figure}

	The improvement of $J_{AYC}$ compared to $S_{AYC}$ with $U_{ref}$ of 6 m/s is 1.56\% lower than that of the WF1 case. On the other hand, the improvement of $J_{AYC}$ compared to $S_{AYC}$ of the Vestas V112 case is 1.90\% higher than that of the WF1 case. In Fig.~\ref{fig:Vestas}(\subref{fig:C_T}), when the wind speed is less than 9 m/s, the $C_T$ of both wind turbines are almost the same. This means the wake characteristics normalized by $U_{ref}$ are also similar. The main reason for the difference in the improvements of $J_{AYC}$ compared to $S_{AYC}$ can be explained using the magnitude of the power gradient. In Fig.~\ref{fig:derivative}(\subref{fig:normalpower}), normalized power curves $\mathcal{P}=P/P_{rated}$ of the Vestas V80 and the Vestas V112 are shown. The denominator is selected as the rated power of each wind turbine for generalization. Fig.~\ref{fig:derivative}(\subref{fig:power_derivative}) shows the gradient of the normalized power curves $\mathcal{P}$.
	
	In Table~\ref{tab:powergradcase}, $\bar{\mathcal{P}}_{up}$ and $\bar{\mathcal{P}}_{down}$ are shown. Subscript ${up}$ means upstream wind turbines where the generated power decreased by AYC and subscript ${down}$ means downstream wind turbines where the generated power is increased by AYC. $\bar{\mathcal{P}}_{up}$ and $\bar{\mathcal{P}}_{down}$ represent normalized power averaged by the number of upstream and downstream wind turbines, respectively. $\bar{\mathcal{P}}_{up}$ is higher than $\bar{\mathcal{P}}_{down}$ because winds through the upstream wind turbines are not much affected by the wake. As shown in Fig.~\ref{fig:derivative}, the power gradients of both wind turbines monotonically increase until 9 m/s. Therefore, upstream wind turbines are in a higher power gradient region than downstream wind turbines. In all cases, $(\overline{\frac{\Delta\mathcal{P}}{\Delta u_{in}}})_{up}$ is higher than $(\overline{\frac{\Delta\mathcal{P}}{\Delta u_{in}}})_{down}$. However, $|\Delta\bar{\mathcal{P}}_{up}|$ is lower than $|\Delta\bar{\mathcal{P}}_{down}|$ because the velocity decrease of an upstream wind turbine which is proportional to $\cos{\gamma}$ is much smaller than the velocity increase of a downstream wind turbine. Deflection of the wake increases as it flows downstream. The difference between $|\Delta\bar{\mathcal{P}}_{up}|$ and $|\Delta\bar{\mathcal{P}}_{down}|$, which is the net power improvement by AYC, is higher when the wind turbines are in a high power gradient region. For the Vestas V80, the power gradient of $U_{ref}=$ 8 m/s case is higher than that of $U_{ref}$= 6 m/s case. Similarly, the power gradient of the Vestas V112 is higher than that of the Vestas V80. As the efficiency of AYC becomes higher, the efficiency of the joint optimization increases.

    \subsection{Number of turbines per area}\label{subsec:Domain size effect}
   
\begin{figure}
    \centering
    \begin{subfigure} {0.49\textwidth}
        \centering
        \includegraphics[width=\textwidth]{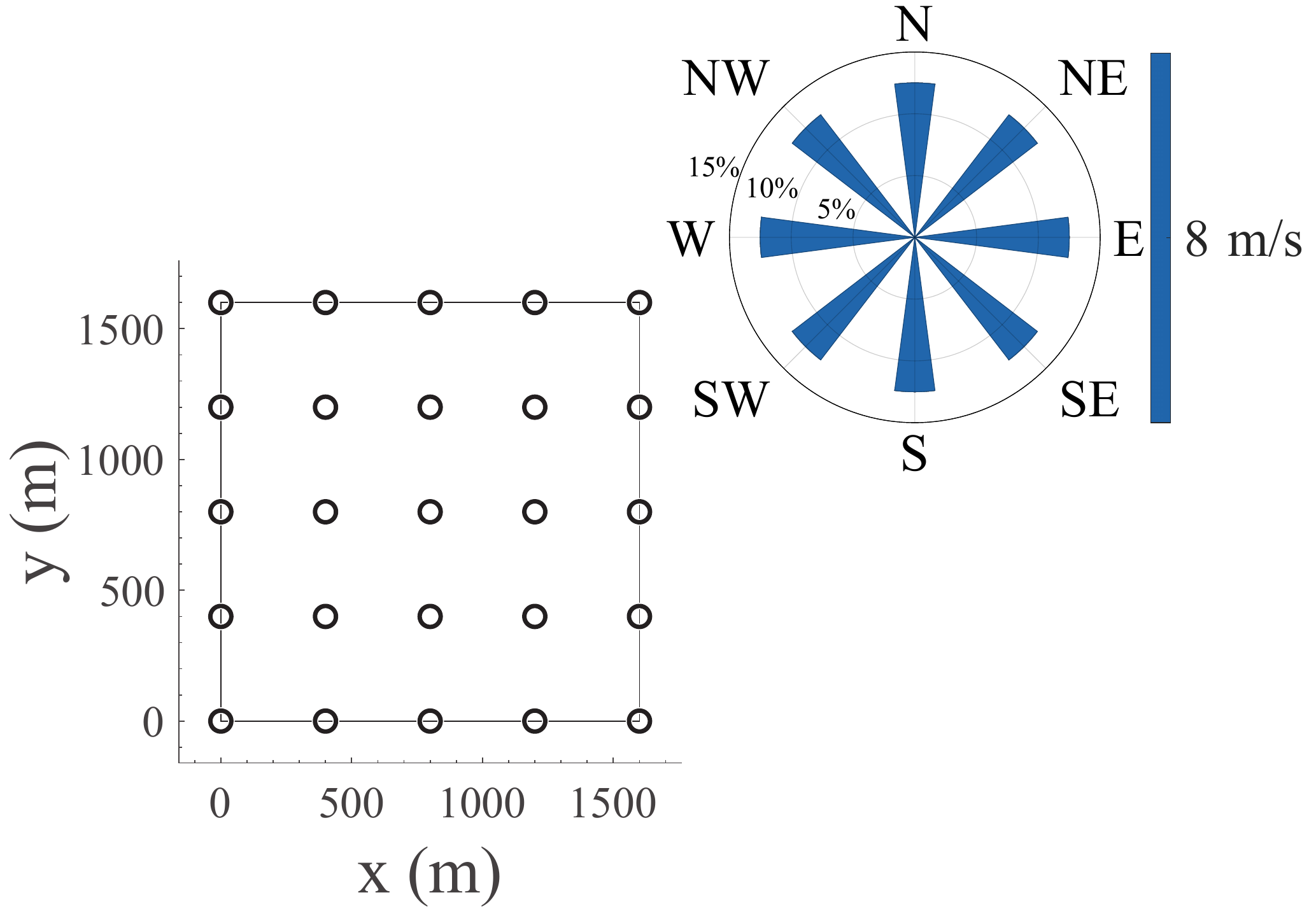}
        \caption{WF1}
        \label{fig:5by5_5D_uniform}
    \end{subfigure}
    \hfill
    \begin{subfigure} {0.49\textwidth}
        \includegraphics[width=\textwidth]{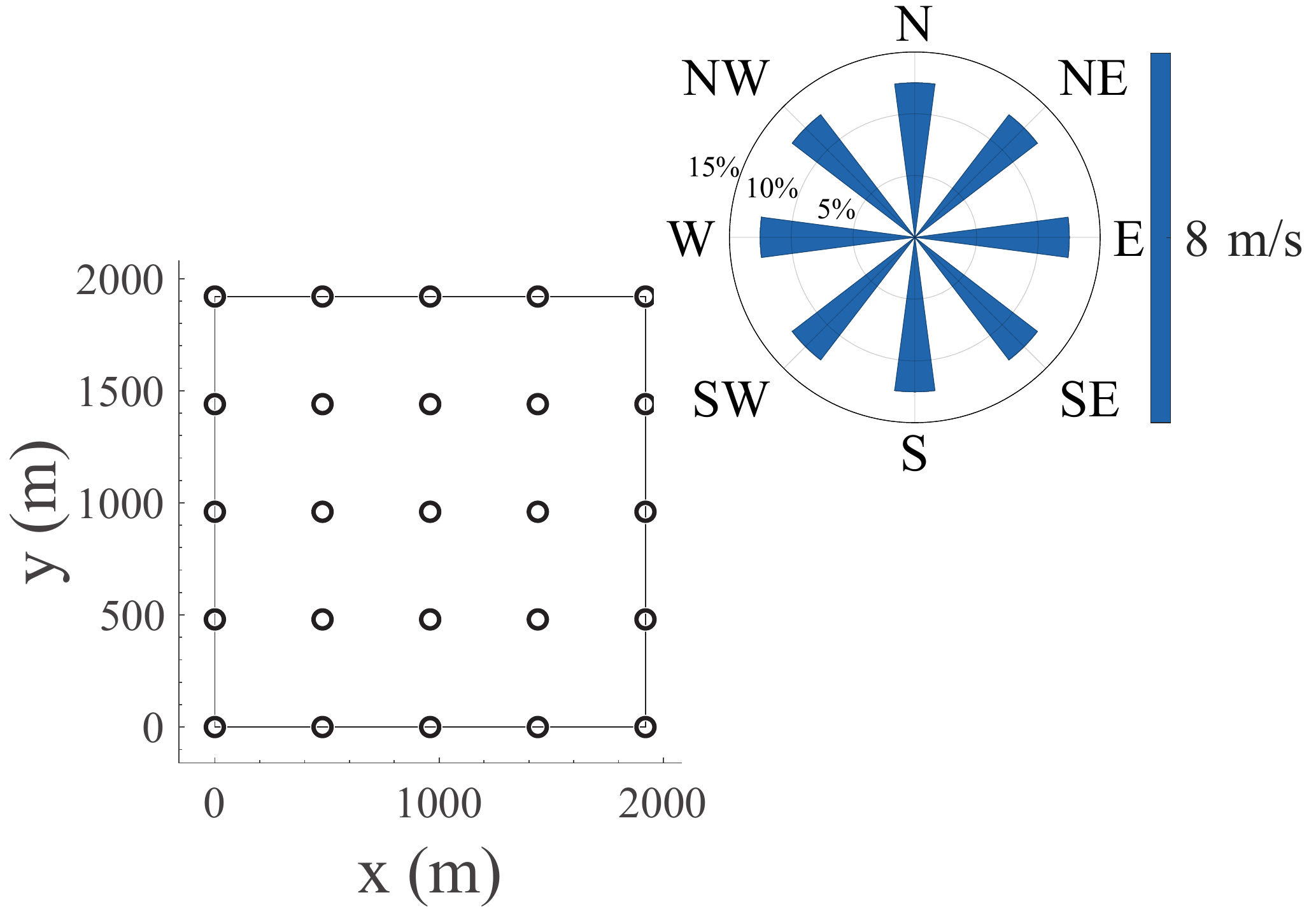}
        \caption{WF2}
        \label{fig:5by5_6D_uniform}
    \end{subfigure}	
    \hfill
    \begin{subfigure} {0.49\textwidth}
        \includegraphics[width=\textwidth]{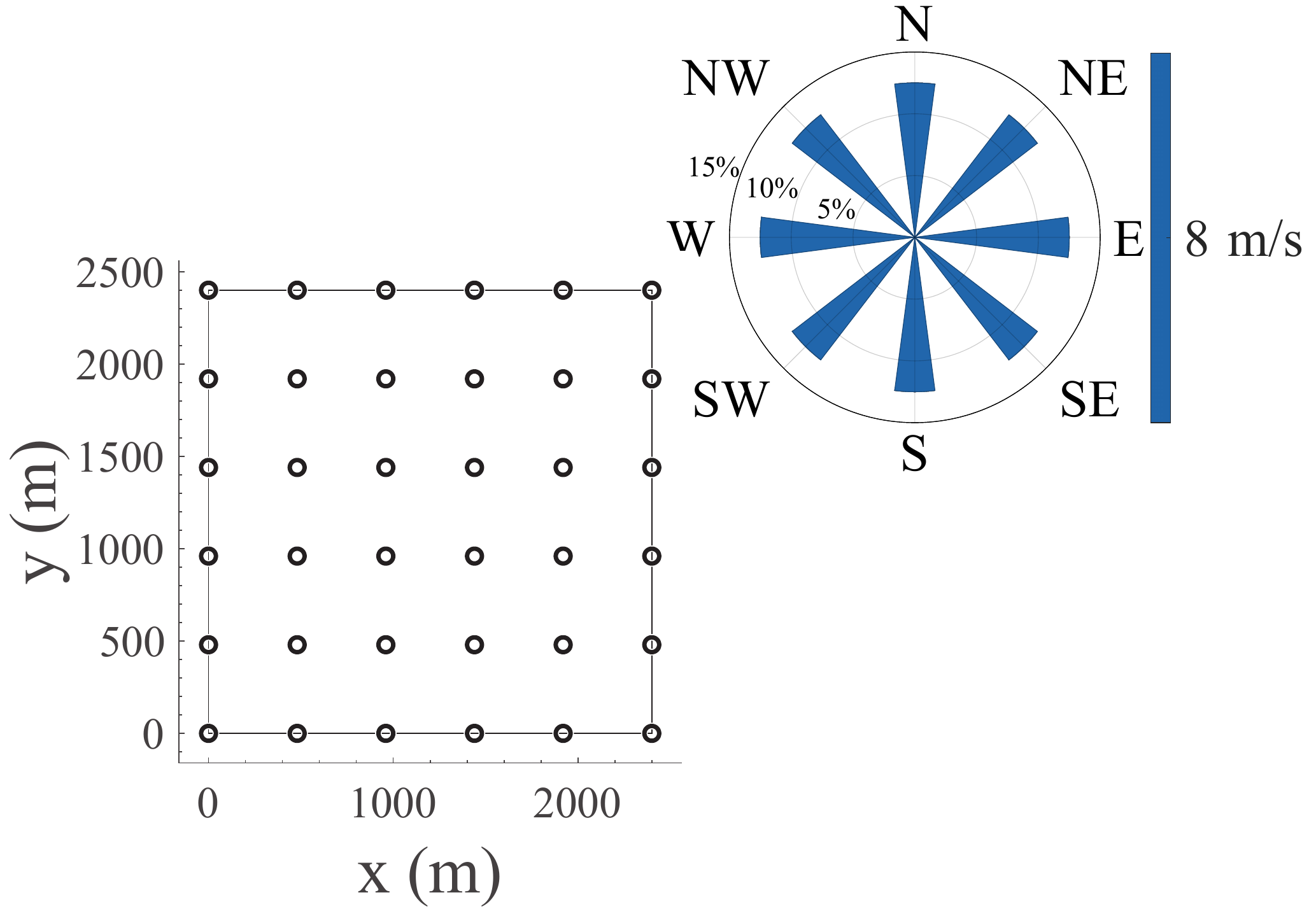}
        \caption{WF3}
        \label{fig:6by6_6D_uniform}
    \end{subfigure}	
    \hfill
    \begin{subfigure} {0.49\textwidth}
        \includegraphics[width=\textwidth]{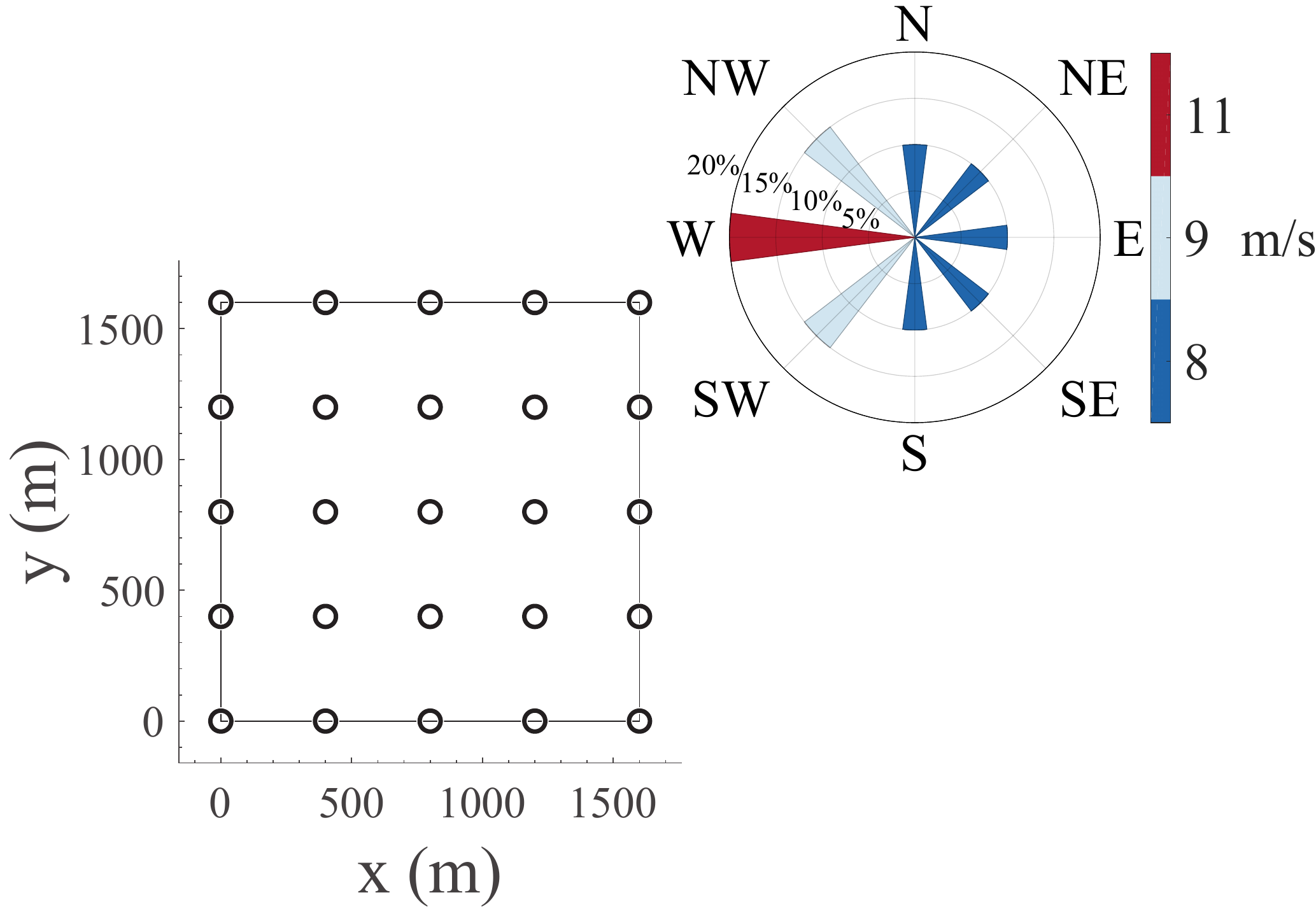}
        \caption{WF4}
        \label{fig:5by5_5D_uneven}
    \end{subfigure}	
    \caption{Four different optimization cases. (a) Optimization of a layout of 25 Vestas V80 wind turbines considering a uniform wind of 8 m/s in a square-shaped wind farm with a side length of 1600 m. The side length is 4$\times$5$D$ where $D$ is the wind turbine diameter (80 m). (b) Optimization of a layout of 25 Vestas V80 wind turbines considering a uniform wind of 8 m/s in a square-shaped wind farm with a side length of 1920 m. The side length is 4$\times$6$D$. (c) Optimization of a layout of 36 Vestas V80 wind turbines considering a uniform wind of 8 m/s in a square-shaped wind farm with a side length of 2400 m. The side length is 5$\times$6$D$. (d) Optimization of a layout of 25 Vestas V80 wind turbines considering uneven winds in a square-shaped wind farm with a side length of 1600 m. The side length is 4$\times$5$D$.}
    \label{fig:WFLO_cases}
\end{figure}

    The grid interval of the checkerboard layout of the WF1 case is 5$D$ as shown in Fig.~\ref{fig:WFLO_cases}(\subref{fig:5by5_5D_uniform}). To investigate the effect of the wind farm area on joint optimization efficiency, the grid interval of the WF2 case is increased to 6$D$, and the number of wind turbines is maintained, as shown in Fig.~\ref{fig:WFLO_cases}(\subref{fig:5by5_6D_uniform}). This leads the average distance between wind turbines of an optimized layout to increase. Table~\ref{tab:Result_WF2WF3WF4} shows that AEP of the WF2 case is improved compared to AEP of the WF1 case because the velocity deficit due to the wake is inversely proportional to the distance. The WF2 case is less affected by the wake than the WF1 case. On the other hand, power improvement using AYC instead of the greedy strategy in the WF2 case is found to be lower than that observed in the WF1 case. As the wake velocity deficit becomes smaller, the maximum power gain of the downstream turbine by AYC decreases, so the AYC efficiency reduces. As a result, improvement of $J_{AYC}$ compared to $S_{AYC}$ of the WF2 case is 1.97\% smaller than that of the WF1 case.
    
    In the WF3 case, the grid interval is the same as in the WF2 case as shown in Fig.~\ref{fig:6by6_6D_uniform}, so the wake deficit of wind turbines is similar. In the WF3 case, the number of wind turbines is increased to 36. Table~\ref{tab:Result_WF2WF3WF4} shows power improvement using AYC instead of the greedy strategy in the WF3 case is lower than that observed in the WF2 case. As the number of wind turbines increases, deflected wake regions due to AYC are more likely to overlap with other downstream wind turbines. Therefore, improving the total AEP of the wind farm using AYC is much more difficult. On the other hand, improvement in $J_{AYC}$ compared to $S_{AYC}$ of the WF3 case is 0.45\% higher than that of the WF2 case. This is because the effect of distributing the wake zone by simply changing the positions of wind turbines is lowered for a large number of wind turbines. The possibility of the downstream wind turbines being in the wake zone, which is parallel to the wind direction, increases. In the sequential optimization problem, the layout is not optimized considering wake deflection, so it cannot efficiently use AYC as in the joint optimization problem.
    
    \subsection{Wind uniformity}\label{subsec:Wind rose effect}
    
    To investigate the effect of wind uniformity on the joint optimization efficiency, uneven winds are supposed in the WF4 case, as shown in Fig.~\ref{fig:WFLO_cases}(\subref{fig:5by5_5D_uneven}). The wind from the west has a magnitude of 11 m/s and the probability of the wind is 20\% and the winds from the northwest and the southwest have a magnitude of 9 m/s and the probability of 15\% in each direction. The winds from the other directions have a magnitude of 8 m/s and the probability of 10\%. In Table~\ref{tab:Result_WF2WF3WF4}, improvement of $J_{AYC}$ compared to $S_{AYC}$ of the WF4 case is 1.44\% lower than that in the WF1 case. This can be explained by analyzing the efficiency of the optimal layout of the joint optimization problem for each wind direction. Table~\ref{tab:Result_WF4} shows the difference in the mean AEP of (1) the west, (2) the northwest and the southwest, and (3) the other directions. Unlike directions (2) and (3), the wind energy of direction (1) is significantly large. This means that especially optimizing the wind farm layout for the wind from the west is advantageous in improving the total AEP of the wind farm. Then, the degrees of freedom of the layout in the other directions is lowered, and the wake effect in those directions is increased. In Table~\ref{tab:Result_WF4}, $J_{AYC}$ of direction (1) accounts for 59.2\% of the sum of mean $J_{AYC}$ of directions (1), (2), and (3). The improvements of AYC compared to the greedy strategy increases in the order of directions (1), (2), and (3) because AEP loss caused by the wake increases in that order. Because improvement of $J_{AYC}$ compared to $S_{AYC}$ of direction (1) in the WF4 case is much less than that in the WF1 case, the efficiency of the joint optimization of the WF4 case is lower than that in the WF1 case.
    
    \section{Concluding remarks}\label{sec:Concluding}
    A particle swarm optimization based optimization method for a wind farm layout with consideration of active yaw control has been developed, and its effectiveness has been analyzed. In the present method, a wind farm layout is optimized with simultaneous consideration of yaw angles for all wind velocities to find a globally optimal layout of the joint optimization problem. An in-house code has been developed to calculate the power production of a wind farm using three-dimensional wake and wake merging models, while the annual energy production of the wind farm has been optimized using locations and yaw angles of wind turbines as control variables. The present method uses an adaptive granularity learning distributed particle swarm optimization algorithm and shows outstanding capability in large-scale joint optimization problems. In the WF1 case, which has been used to optimize the layout of 25 Vestas V80 wind turbines considering a uniform wind velocity of 8 m/s in a square-shaped wind farm with a side length of 1600 m, the annual energy production of the present method has been 1.62\% higher than that of the previous method. 
    
    Efficiencies of joint optimization problems of the following conditions have been analyzed. Firstly, the effect of the power curve gradient on optimization has been analyzed. Two modified WF1 cases which have different free stream velocities and the wind turbine model have been optimized. The larger the power curve gradient, the higher the efficiency of active yaw control. When the Vestas V112 was used for the WF1 case, $J_{AYC}$ has been increased by 4.32\% than $S_{AYC}$. Secondly, the WF2 case and the WF3 case, which have different wind farm areas and the number of wind turbines, have been optimized. It has been shown that the efficiency of joint optimization increases as the size of the wind farm decreases and as the number of wind turbines increases. Finally, the effect of the uniformity of winds has been analyzed. In the WF4 case, which supposed uneven winds, because the layout has been optimized to the direction where the wind energy is rich, the efficiency of the joint optimization has been lower than in the WF1 case.
    
    The annual energy production of the jointly optimized layout using the proposed method has been higher than that of the sequentially optimized layout. When designing a jointly optimized layout, the annual energy production with simultaneous consideration of all wind velocities has been higher than that of separately considering wind velocities. The presented method can adopt massive variables and constraints of layout optimization problems and can be widely applied to improve actual wind farm design.
    
	\section*{Acknowledgements}
	The work was supported by the National Research Foundation of Korea (NRF) under the Grant Number NRF-2021R1A2C2092146 and the Korea Electric Power Corporation (KEPCO) under the Grant No. R20XO02-21.

\newpage
\bibliographystyle{elsarticle-num}
\biboptions{sort&compress}
\bibliography{PSO_with_AYC}

\newpage
\pagebreak
\clearpage
\listoftables

\pagebreak
\clearpage
\begin{table}
\caption{Results of the WF1 case optimizing the layout of 25 Vestas V80 wind turbines considering uniform wind with a velocity of 8 m/s in a square-shaped wind farm with a side length of 1600 m. Configurations of the WF1 are shown in Fig.~\ref{fig:WFLO_cases}(\subref{fig:5by5_5D_uniform}). $S$ represents the annual energy production (AEP) of the sequentially optimized result, and $J$ represents AEP of the jointly optimized result. Subscript $\Theta$ means turbines are aligned perpendicular to the wind direction. Subscript $AYC$ is the active yaw control (AYC) method where the yaw angles of turbines are optimized to maximize AEP of the wind farm. Superscripts $DBHM$ and $AGLD$ represent that the results are optimized using the decomposition-based hybrid method (DBHM) and the adaptive granularity learning distributed particle swarm optimization (AGLDPSO) algorithms, respectively.}
\centering 
\begin{tabular}{ccccccc} 
\hline 
 & $S_{\Theta}$ & $S_{AYC}$ & $J_{\Theta}^{DBHM}$ & $J_{AYC}^{DBHM}$ & $J_{\Theta}^{AGLD}$ & $J_{AYC}^{AGLD}$ \\ [1ex]
\hline 
AEP (GWh) & 168.46 & 170.61 & 167.19 & 171.96 & 167.98 & 174.74 \\[1ex]
Improvement (\%) & - & 1.28 & -0.75 & 2.08 & -0.29 & 3.73 \\[1ex]
\hline 
\end{tabular}
\label{tab:5by5_5D}
\end{table}

\pagebreak
\clearpage
\begin{table}
\caption{Results of the modified WF1 cases. Both cases optimize the layout of 25 wind turbines in a square-shaped wind farm with a side length of 1600 m. $U_{ref}$ of the case (a) is $6$m/s and turbines of the case (a) are Vestas V80. $U_{ref}$ of the case (b) is $8$m/s and the turbines of the case (b) are Vestas V112.}
\centering 
\begin{tabular}{ccccccc} 
\hline 
 & \multicolumn{2}{c}{Modified Variables} & $S_{\Theta}$ & $S_{AYC}$ & $J_{\Theta}$ & $J_{AYC}$ \\ [1ex]
\hline 
& & AEP (GWh) & 69.35 & 70.62 & 68.46 & 71.23 \\[-1ex]
\raisebox{1.5ex}{(a)} & \raisebox{1.5ex}{$U_{ref}=6$m/s} & Improvement (\%) & - & 1.83 & -1.29 & 2.70 \\[1ex]
\hline 
& & AEP (GWh) & 323.20 & 333.72 & 321.04 & 348.14 \\[-1ex]
\raisebox{1.5ex}{(b)} & \raisebox{1.5ex}{Vestas V112} & Improvement (\%) & - & 3.26 & -0.67 & 7.72 \\[1ex]
\hline
\end{tabular}
\label{tab:modifiedWF1}
\end{table}

\pagebreak
\clearpage
\begin{table}
\caption{Results of $J_{AYC}$ of the WF1 case and of the modified WF1 cases, which have different $U_{ref}$ and the wind turbine model, respectively. $\mathcal{P}$ is the power of a wind turbine normalized by its rated power. Subscript ${up}$ means upstream wind turbines where the generated power decreases due to active yaw control (AYC) and subscript ${down}$ means downstream wind turbines where the generated power increases due to AYC. $\bar{\mathcal{P}}_{up}$ and $\bar{\mathcal{P}}_{down}$ are normalized power averaged by the number of upstream and downstream wind turbines, respectively. $\Delta$ means the amount of value changed by AYC. $(\overline{\frac{\Delta\mathcal{P}}{\Delta u_{in}}})$ represents the averaged derivative of power.}
\centering 
\begin{tabular}{ccccccc} 
\hline 
 & $\bar{\mathcal{P}}_{up}$ & $\Delta\bar{\mathcal{P}}_{up}$ & $(\overline{\frac{\Delta\mathcal{P}}{\Delta u_{in}}})_{up}$ & $\bar{\mathcal{P}}_{down}$ & $\Delta\bar{\mathcal{P}}_{down}$ & $(\overline{\frac{\Delta\mathcal{P}}{\Delta u_{in}}})_{down}$ \\ [1ex]
\hline 
WF1 & 0.438 & -0.020 & 0.152 & 0.330 & 0.055 & 0.135\\[1ex]
\hline 
$U_{ref}=6 m/s$ & 0.178 & -0.009 & 0.086 & 0.134 & 0.022 & 0.077\\[1ex]
\hline 
Vestas V112 & 0.615 & -0.059 & 0.192 & 0.378 & 0.133 & 0.167\\[1ex]
\hline
\end{tabular}
\label{tab:powergradcase}
\end{table}

\pagebreak
\clearpage
\begin{table}
\caption{Results of the WF2, WF3, and WF4 cases. For the WF2 case, the layout of 25 Vestas V80 wind turbines is optimized by considering a uniform wind of 8 m/s in a square-shaped wind farm with a side length of 1920 m. For the WF3 case, the layout of 36 Vestas V80 wind turbines is optimized by considering a uniform wind of 8 m/s in a square-shaped wind farm with a side length of 2400 m. For the WF4 case, the layout of 25 Vestas V80 wind turbines is optimized by considering uneven winds in a square-shaped wind farm with a side length of 1600 m. The configurations of the WF2, WF3, and WF4 are shown in Fig.~\ref{fig:WFLO_cases}(\subref{fig:5by5_6D_uniform}), Fig.~\ref{fig:WFLO_cases}(\subref{fig:6by6_6D_uniform}), and Fig.~\ref{fig:WFLO_cases}(\subref{fig:5by5_5D_uneven}), respectively.}
\centering 
\begin{tabular}{cccccc} 
\hline 
& & $S_{\Theta}$ & $S_{AYC}$ & $J_{\Theta}$ & $J_{AYC}$ \\ [1ex]
\hline 
& AEP (GWh) & 180.44 & 182.53 & 178.15 & 183.35 \\[-1ex]
\raisebox{1.5ex}{WF2} &Improvement (\%) & - & 1.16 & -1.27 & 1.61 \\[1ex]
\hline
& AEP (GWh) & 249.79 & 251.93 & 248.97 & 254.21 \\[-1ex]
\raisebox{1.5ex}{WF3} &Improvement (\%) & - & 0.86 & -0.33 & 1.77 \\[1ex]
\hline
& AEP (GWh) & 233.96 & 236.48 & 232.06 & 238.80 \\[-1ex]
\raisebox{1.5ex}{WF4} &Improvement (\%) & - & 1.08 & -0.81 & 2.07 \\[1ex]
\hline
\end{tabular}
\label{tab:Result_WF2WF3WF4}
\end{table}

\pagebreak
\clearpage
\begin{table}

\caption{AEP of the WF4 case is divided by directions as follows: (1) AEP when the wind direction is west, (2) mean AEP when the wind directions are northwest and southwest, and (3) mean AEP when the wind comes from the other directions.}
\centering 
\begin{tabular}{ccccccc} 
\hline 
Wind & & $S_{\Theta}$ & $S_{AYC}$ & $J_{\Theta}$ & $J_{AYC}$ & $J_{AYC}/({\sum{J_{AYC}}})$ \\ [1ex]
\hline 
& AEP (GWh) & 77.87 & 78.09 & 76.73 & 78.41 & \\[-1ex]
\raisebox{1.5ex}{(1)} &Improvement (\%) & - & 0.29 & -1.46 & 0.70 & \raisebox{1.5ex}{59.2\%}
\\[1ex]
\hline
& AEP (GWh) & 36.26 & 36.70 & 35.62 & 36.76 & \\[-1ex]
\raisebox{1.5ex}{(2)} &Improvement (\%) & - & 1.21 & -1.77 & 1.38 & \raisebox{1.5ex}{27.7\%}
\\[1ex]
\hline
& AEP (GWh) & 16.71 & 17.00 & 16.82 & 17.37 & \\[-1ex]
\raisebox{1.5ex}{(3)} &Improvement (\%) & - & 1.71 & 0.62 & 3.95 & \raisebox{1.5ex}{13.1\%} \\[1ex]
\hline
\end{tabular}
\label{tab:Result_WF4}
\end{table}

\newpage
\pagebreak
\clearpage
\listoffigures

\pagebreak
\clearpage

\pagebreak
\clearpage

\pagebreak
\clearpage

\pagebreak
\clearpage

\pagebreak
\clearpage

\pagebreak
\clearpage

\pagebreak
\clearpage

\pagebreak
\clearpage

\end{document}